\documentclass[11pt]{article}
\usepackage{amstext,amssymb,amsmath,amsbsy}

\usepackage{hyperref}
\usepackage{amsfonts}
\usepackage{verbatim}
\usepackage{amsmath}
\usepackage{amsthm}
\usepackage{enumerate}
\usepackage{graphicx}
\usepackage{color}

\newcommand{\sg}{\sigma}

\newcommand{\ga}{\hat \eps}

\newcommand{\Og}{\Omega}

\newcommand{\pdh}{\partial}

\newcommand{\Ga}{\Gamma}

\newcommand{\di}{\mbox{div}}

\newcommand{\mG}{\mathcal{G}}

\newcommand{\mD}{\mathcal{D}}

\newcommand{\me}{{\bf e}}

\newcommand{\mw}{{\bf d}}

\newcommand{\mH}{\mathcal{H}}

\newcommand{\mM}{\mathcal{M}}

\newcommand{\ag}{\alpha}

\newcommand{\pn}{\partial_{n}}

\newcommand{\intl}{\int\limits}

\newcommand{\mF}{\mathcal{F}}

\newcommand{\mf}{k}

\newcommand{\ba}{\begin{eqnarray}}

\newcommand{\bas}{\begin{eqnarray*}}

\newcommand{\ea}{\end{eqnarray}}

\newcommand{\eas}{\end{eqnarray*}}

\newcommand{\rN}{\mathbb{R}}

\newcommand{\vn}{n}

\newcommand{\mq}{q}

\newcommand{\ms}{s}

\newtheorem{theorem}{Theorem}
\newtheorem{lemma}{Lemma}

\newtheorem{proposition}{Proposition}
\newtheorem{corollary}{Corollary}
\newtheorem{remark}{Remark}

\newtheorem{definition}{Definition}
\numberwithin{equation}{section}

\newcommand{\proofend}{\hfill $\Box$ }

\newcommand{\dsp}{\displaystyle}

\newcommand{\supp}{\operatorname{supp}}

\newcommand{\dive}{\operatorname{div}}

\newcommand{\eps}{\varepsilon}

\newcommand{\loc}{_{loc}}

\newcommand{\mN}{\mathbb{N}}

\newcommand{\mR}{\mathbb{R}}

\newcommand{\mA}{\mathcal{A}}

\newcommand{\mc}{\mathrm{c}}

\title{Generalized Impedance Boundary Conditions for Strongly Absorbing Obstacles: the full Wave Equations}
\author{Hoai-Minh Nguyen \footnote{EPFL SB MATHAA CAMA, Station 8,  CH-1015 Lausanne, Switzerland, E-mail: \texttt{hoai-minh.nguyen@epfl.ch} and School of Mathematics, University of Minnesota, MN, 55455, E-mail: \texttt{hmnguyen@math.umn.edu}. }
~  and Linh V.~Nguyen\footnote{Department of Mathematics, University of Idaho,
Moscow, Idaho 83843, USA, \texttt{lnguyen@uidaho.edu}. }}

\date{\today}

\begin{document}
\maketitle

\begin{abstract}  This paper is devoted to the study of  the generalized impedance boundary conditions (GIBCs) for a strongly absorbing obstacle in the {\bf time} regime in two and three dimensions.  The GIBCs in the time domain are heuristically derived from the corresponding conditions in the time harmonic regime. The latters are frequency dependent except the one of order 0; hence the formers are  non-local in time in general. The error estimates in the time regime can be derived  from  the ones in the time harmonic regime when the frequency dependence is well-controlled. This idea is originally due to Nguyen and Vogelius in \cite{NguyenVogelius2} for the cloaking context. In this paper, we present the  analysis to the GIBCs of orders 0 and 1. To implement the ideas in \cite{NguyenVogelius2}, we revise and extend the work of Haddar, Joly, and Nguyen in \cite{HJNg1}, where the GIBCs were investigated for a fixed frequency in three dimensions. Even though we heavily follow the strategy in \cite{NguyenVogelius2}, our analysis on the stability contains new ingredients and ideas. First, instead of considering the difference between solutions of the exact model and the approximate model, we consider the difference between their derivatives in time. This simple idea helps us to avoid the machinery used in \cite{NguyenVogelius2}  concerning the integrability with respect to frequency in the low frequency regime. Second, in the high frequency regime, the Morawetz multiplier technique used in \cite{NguyenVogelius2} does not fit directly in our setting. Our proof makes use of a result by H\"ormander in \cite{Hor}. Another important part of the analysis in this paper is the well-posedness in the time domain for the approximate problems imposed with GIBCs on the boundary of the obstacle, which are non-local in time. 

\end{abstract}

\tableofcontents

\section{Introduction and statement of the main results}\label{S:Intro}

The computation of electromagnetic scattering from an arbitrary obstacle has been an active research area for many decades. One technique is to replace the exact model inside the obstacle by appropriate  boundary conditions on its surface; hence the problem of determining the external electromagnetic fields can be solved without considering the fields inside the obstacle (see, e.g., \cite{HoppeSami,SeniorVolakis}). These boundary conditions are called {\it Generalized Impedance Boundary Conditions} (GIBCs). The first GIBC for a highly absorbing obstacle (highly conducting obstacle) was proposed by Leontovich (see, e.g., \cite{Leon}) and was extended later by  Rytov in \cite{Rytov}. Antoine, Barucq, and Vernhet in  \cite{AntoineBarucqVernhet}, using the technique of pseudo-differential equations (following the ideas of Engquist and Majda in \cite{EngquistMajda}), implemented a new derivation of such conditions. Recently, Haddar, Joly, and Nguyen in \cite{HJNg1}  revisited these GIBCs for the Helmholtz equation. More precisely,  the authors first proposed a new construction of GIBCs  which is based on an ansatz for the asymptotic expansion of exact solutions. They then developed mathematical tools, based on compactness arguments, to establish error estimates up to order 3. Related works are the edgy current problem studied by  MacCamy and Stephan in \cite{MacSte} and the study of  the GIBCs for highly conducting obstacle for the Maxwell system in the time harmonic regime by Haddar, Joly, and Nguyen in \cite{HJNg2} (see also the work of Caloz, Dauge, Faou, and P\'eron in \cite{CarlosPeron1}) and references therein. 

\medskip
The study of GIBCs for highly conducting obstacle has been though studied extensively in the literature, the rigorous study of GIBCs for highly absorbing obstacle in the time regime is not known to our knowledge. The lack of the study in the time regime is not special for this context but a common problem in the study of acoustic and electromagnetic waves since problems in the time regime involve the interaction of all frequency and hence they are harder to analyze. 

\medskip
The goal of this  paper is to provide the analysis of the GIBCs for highly conducting obstacle in the {\bf time regime} in two and three dimensions. Heuristically, these are obtained by taking the inverse Fourier transform of the corresponding conditions in the time harmonic regime with respect to frequency. Since the GIBCs in the time harmonic regime are frequency dependent, the ones in the {\bf time} regime are {\bf non-local} with respect to time.  
The error estimates in the time regime can be derived  from  the ones in the time harmonic regime, when the frequency dependence is well-controlled. This idea is originally due to Nguyen and Vogelius in \cite{NguyenVogelius2} used for the cloaking context. To implement this idea, we revise and extend the work of Haddar, Joly, and Nguyen in \cite{HJNg1}, where the GIBCs for the time harmonic regime  were investigated for a fixed frequency in three dimensions. Even though, we follow the strategy in \cite{NguyenVogelius2}, our analysis on the stability contains new ingredients and ideas. First, instead of considering the difference between solutions of the exact model and the approximate model in the time harmonic regime, we consider the difference between their derivative in time. This simple idea helps us to avoid the machinery used in \cite{NguyenVogelius2}  concerning the integrability with respect to frequency in the low frequency regime. The proof of the stability in the low frequency regime is based on a compactness argument as in \cite{HJNg1} and uses ideas in \cite{Ng-cl-1}. Second, the compactness argument is not appropriate in the high frequency regime; moreover, the well-known Morawetz multiplier technique used in \cite{NguyenVogelius2} does not fit directly in our setting,  mainly due to the structure of the GIBCs. This new technical challenge distinguishes our work from \cite{NguyenVogelius2}. To tackle it, we use essentially a result by H\"ormander in \cite{Hor} (Lemma~\ref{lem-Hor}). Another important part of the proof is the well-posedness in the time domain of the approximate problems imposed with GIBCs on the boundary of the obstacle, which are non-local in time.   

\medskip

In this paper, we concentrate only on the GIBCs of order {\bf 0 and 1}. Even though our approach also works for the GIBCs of orders 2 and 3, we are still not able to obtain the optimal expected estimates as in the time harmonic regime due to the complexity of the structures of the GIBCs in these cases. We postpone the study of these conditions to our future work \cite{NgNg-1}.

\medskip
There are many other interesting situations in which the asymptotic expansions have been investigated in the time harmonic setting. For example, the thin coating effect in e.g., \cite{BenLem}, the wave propagation in media with thin slots  in e.g., \cite{JolTor}, the wave propagation across thin periodic interfaces in e.g., \cite{DelHarJol} and references therein.  We hope that our analysis also sheds light to the situations mentioned.

\medskip


We next describe the problem in more details. Let $\Og$ be a smooth bounded domain in $\rN^d$ ($d=2,3$) with boundary $\Gamma : = \partial \Omega$ and
let $f \in L^{2}(\rN_+ \times \rN^d)$ be such that $\supp f \subset [0,T] \times (B_{R_0} \setminus \overline{\Og}) $ for some fixed $R_0>0$ and $T>0$. 
Here and in what follows, $B_r$ denotes the ball of radius $r$ centered at the origin. Let  $u^{\eps} \in L_{\loc}^\infty\big([0, \infty); H^{1}(\rN^d)\big)$ with $\partial_{t} u^{\eps} \in L_{\loc}^\infty\big([0, \infty); L^{2}(\rN^d)\big)$ be the unique weak solution to the problem
\begin{equation}\label{E:wave} 
\left\{\begin{array}{ll} \dsp \pdh_{tt} u^{\eps}(t,x) - \Delta u^{\eps}(t,x) +\sigma_\eps(x) \; \pdh_t u^{\eps}(t,x)  = f(t,x),&\hskip 0 pt (t,x) \in \rN_+ \times \rN^d, \\[6 pt] u^{\eps}(0,x)=\pdh_t u^{\eps}(0,x) =0,&\hskip 0 pt x \in \rN^d,  \end{array}\right. 
\end{equation}
where 
\begin{equation*}
\sg_\eps(x) = \left\{\begin{array}{cl} 0 & \mbox{ if } x \in \rN^{d} \setminus \Og, \\[6 pt] \dsp \frac{1}{\eps^2}& \mbox{ if } x\in \Og. \end{array} \right.
\end{equation*}
for some $\eps > 0$ small.  Roughly speaking, the absorption of $\Omega$ is of order $1/ \eps^{2}$.  We consider here the case in which the initial conditions are zero. The general case could be treated similarly as discussed in the cloaking setting in \cite{NguyenVogelius2}.  

\medskip
Let $\hat u_\eps(k, x)$ and $\hat f(k, x)$ be the Fourier transform of $u_\eps$  and $f$ with respect to time respectively, i.e., \footnote{We extend these function by 0 for $t < 0$.}
\begin{equation}\label{def-F-u}
\hat u_\eps(k, x) = {\cal F} (u_\eps) (k, x) = \frac{1}{\sqrt{2 \pi}} \int_{\rN} u_\eps(t, x) \, e^{i\, k \, t} \, dt
\end{equation}
and 
\begin{equation}\label{def-F-f}
\hat f(k, x) = {\cal F} (f) (k, x) = \frac{1}{\sqrt{2 \pi}} \int_{\rN} f(t, x) \, e^{i\, k \, t} \, dt.
\end{equation}
Then for almost every $k > 0$, $\hat u_\eps(k, x) \in H^1_{\loc}(\mR^d)$ be the unique solution to the equation 
\begin{equation}\label{eq-hat-u}
\Delta \hat u_\eps(k, x)   + \mf^2 \hat u_\eps(k, x)  +i \, \mf \, \sg_\eps \; \hat u_\eps(k, x)  = \hat f(k, x) ,\quad \quad \mbox{ in } \rN^d, 
\end{equation}
which satisfies the outgoing condition 
\begin{equation}\label{OGC}
\pdh_r u_\eps - i \mf  \; u_\eps   = o \big(r^{- \frac{(d-1)}{2}}\big), \quad \mbox{ as } r =|x| \to \infty. 
\end{equation}
This fact is formulated later in Lemma \ref{lem-outgoing-3} whose proof has root in \cite[Theorem A.1]{NguyenVogelius2}. 

\medskip
Set \footnote{$\alpha^2 = -i$.} 
\begin{equation} \label{notation}
\ga= \frac{\eps}{\sqrt{\mf}},  \quad   \quad \alpha =  \frac{\sqrt{2}}{2} - i \frac{\sqrt{2}}{2}, 
\end{equation}
and 
\begin{equation} \label{E:BFreq} 
\mD_\ell^{\ga} = \left\{\begin{array}{ll} 0 & \mbox{ for }  \ell =0, \\[6 pt] 
\dsp \frac{\ga}{\alpha} &\mbox{ for } \ell=1.\end{array}\right.
\end{equation}
It is proved in  \cite{HJNg1} that the GIBCs of order 0 and 1 corresponding to  \eqref{eq-hat-u} are 
\begin{eqnarray}  \label{E:GIBC-frequency}
v + \mD_{\ell}^{\ga} \, \pn v = 0, &\hskip 0 pt  \mbox{ on } \Ga. 
\end{eqnarray}  
Here and in what follows, $n=n(x)$ is the unit normal vector {\bf directed into} $\Omega$ on $\Ga$. More precisely, let $s \in H^1(\mR^d)$ with support in $B_{R_0} \setminus \Omega$ and  let $v^\eps \in H_{\loc}^1(\mR^d)$ and $v_{\ell}^{a} \in H^1_{\loc}(\rN^d \setminus \Omega)$ be the unique outgoing solutions to the problems
\begin{equation}\label{def-v} \Delta v^{\eps}  + \mf^2 v^{\eps} +i \, \mf \, \sg_\eps \; v^{\eps} = \ms,\quad \quad \mbox{ in } \rN^d, \end{equation}
and
\begin{equation}\label{def-va}\left\{\begin{array}{ll} \Delta v_{\ell}^{a} +\mf^2 \; v_{\ell}^{a} = \ms,&\hskip 0 pt  \mbox{ in } \rN^d \setminus \Omega,\\[6 pt]
v_{\ell}^{a} + \mD_{\ell}^{\ga}  \; \pn v_{\ell}^{a} =0 & \mbox{ on } \Gamma.\end{array}
\right.
\end{equation} 
Haddar, Joly, and Nguyen \cite[Theorem 3]{HJNg1} proved that, for any $R>0$, \footnote{In \cite{HJNg1}, the authors considered the bounded setting. However, their method also implies the results mentioned here.} 
\begin{equation}\label{known-C}
\|v^{\eps} - v_{\ell}^{a}\|_{H^1(B_R \setminus \Og)} \leq C(\mf, R) \; \eps^{\ell+1} \| s\|_{H^m(\rN^d)}, 
\end{equation}
for some positive constant $C(\mf, R)$ and for some $m>0$ large enough. The dependence on $\mf$ of $C(\mf, R)$ in \cite{HJNg1} is not explicit. 

\medskip
We are now ready to {\bf heuristically} derive the GIBCs for \eqref{E:wave} by  taking the inverse Fourier transform of the GIBCs in the time harmonic regime with respect to frequency.  We have

\medskip
\noindent {\bf GIBC of order 0}:
\begin{equation} \label{E:D} G^\eps_{0}(v) := v=0, \quad \mbox{ on } \rN_+ \times \Ga.\end{equation}
This is clear from  \eqref{E:BFreq} and \eqref{E:GIBC-frequency} with $\ell =0$. 

\medskip
\noindent {\bf GIBC of order 1}: \begin{equation} \label{E:GIBC-wave}
G^\eps_{1}(v) :=  \pn v + B_{1}^{\eps} \, v = 0   \mbox{ on } \rN_+ \times \Ga,
\end{equation}
where 
\begin{equation}\label{E:B}
\big(B^\eps_1 \, v \big) (t,x) : =  \frac{1}{\sqrt{\pi} \, \eps} \, \int\limits_0^t \frac{ \pdh_t v(\tau,x)}{\sqrt{t-\tau}} \, d\tau.
\end{equation}
The derivation goes as follows. For $\ell=1$, condition \eqref{E:GIBC-frequency} reads as:
\begin{eqnarray*}
\pn v + \frac{\ag \, \sqrt{\mf}}{\eps} \; v= 0, &\hskip 0 pt  \mbox{ on } \Ga,
\end{eqnarray*}
or equivalently 
\begin{eqnarray*}
\pn v + \frac{1}{\eps} \frac{1}{\ag \, \sqrt{\mf}} \; (-i k) v  = 0, &\hskip 0 pt  \mbox{ on } \Ga.
\end{eqnarray*}
The condition \eqref{E:GIBC-wave} is now a consequence of the fact (see,  e.g.,  \cite[p. 171]{GelfandShilov})
\begin{equation}\label{E:conv}\mF (\varphi)(\mf) = \frac{1}{\sqrt{2 \pi}} \; \frac{\sqrt{\pi}}{\ag \sqrt{k}}, \quad \mbox{ where } \quad  
  \varphi(t)= \left\{\begin{array}{cl} \dsp \frac{1}{\sqrt{t}} &  \mbox{ if } t>0, \\[6 pt]  0 & \mbox{ if } t \leq 0. \end{array}  \right.\end{equation}

We have heuristically derived the GIBCs  of orders $0$ \eqref{E:D} and $1$ \eqref{E:GIBC-wave}  for the full wave equation \eqref{E:wave}. Similarly, one can obtain the GIBCs of orders $2,3$ for \eqref{E:wave} from the corresponding ones in the time harmonic regime obtained  in \cite{HJNg1}. However, such conditions are more complicated. We have not been able yet to obtain  the optimal  expected estimates for them as in the ones in the time harmonic regime. We postpone the study of these conditions to our future work \cite{NgNg-1}.


\medskip 
The goal of this paper is to establish error estimates for \eqref{E:D} and \eqref{E:GIBC-wave}. More precisely, we prove
\begin{theorem}\label{T:Main}
Let $d=2, 3$,  $\ell=0,1$, $T>0$, $R_0>0$, and let $G^\eps_{\ell}$ be defined in \eqref{E:D} and \eqref{E:GIBC-wave}. Assume that $f \in L^2\big([0, \infty) \times \rN^{d} \big)$ with $\supp f \subset [0, T] \times (B_{R_0} \setminus \overline \Og)$. 
There exists a unique weak solution $u_{\ell}^{a} \in L_{\loc}^\infty([0,\infty); H^1(\rN^d \setminus \Og))$ with $\partial_t u_{\ell}^{a} \in L_{\loc}^\infty([0,\infty); L^2(\rN^d \setminus \Og))$ to 
\begin{equation}\label{E:wave-ua}
\left\{\begin{array}{ll}
\partial^2_{tt} u_\ell^a - \Delta u_\ell^a = f & \mbox{ in }\rN_+ \times \mR^d \setminus \overline {\Omega}, \\[6pt]
G^\eps_\ell (u_\ell^a) = 0 & \mbox{ on } \mR_+ \times \Gamma, \\[6pt]
\partial_t u_\ell^a (0, \cdot) = u_\ell^a (0, \cdot) = 0 & \mbox{ in } \mR^d \setminus  \overline {\Omega};
\end{array} \right.
\end{equation}
moreover,
\begin{equation}\label{est-stability}
\int\limits_{\rN^{d} \setminus \Omega} |\nabla u_{\ell}^{a}(t, x)|^{2} + |\pdh_t u_{\ell}^{a}(t, x)|^{2} \, dx 
 \leq C\, t \, \|f\|_{L^{2}([0,t] \times \rN^d)}^{2} \quad \forall \, t  \ge 0. 
\end{equation}
Assume in addition that $\Omega$ is {\bf star-shaped} and $f \in C^\infty\big((0, \infty) \times \rN^{d} \big)$ with $\supp f \subset \subset (0, + \infty) \times (B_{R_0} \setminus \overline \Og)$. Then,  for any $t>0$ and $K \subset \subset \rN^d \setminus \overline \Omega$ \footnote{Roughly speaking, $K$ is bounded and away from the boundary of the obstacle.}, there is a positive constant $C$ independent of $\eps$ and $f$,  such that, for some integer $m$ \footnote{The integer $m$ can be chosen as follows $m=13$ if $\ell =0$ and $m=16$ if $\ell =1$. Assume in addition that $\supp f \cap \bar \Omega = \O$. Then $m$ can be chosen as follows $m=8$ if $\ell =0$ and $m=9$ if $\ell=1$; however the constant $C$ in \eqref{est-stability} now depends on the distance between $\supp f$ and $\bar \Omega$ (see Footnote \ref{footnote1}).}. 
\begin{equation}\label{k-wave}
\|u^{\eps} -u_{\ell}^{a}\|_{L^\infty \big([0,  t]; H^{1}(K)\big)}  \leq C \; \eps^{\ell+1} \|f\|_{H^m(\rN_+ \times \rN^d)}. 
\end{equation}
\end{theorem}

The following definition of the weak solutions, which is motivated from the standard concept of weak solutions, is used in Theorem~\ref{T:Main}.
\begin{definition} Let $d=2, 3$ and $\ell=0, 1$. We say a function 
\begin{equation*}
u_{\ell}^{a} \in L_{\loc} ^\infty([0, \infty); H^1(\rN^d \setminus \Omega)) \mbox{ with } \partial_t u_{\ell}^{a} \in L_{\loc}^\infty([0, \infty); L^2(\rN^d \setminus \Omega))
\end{equation*}
is a weak solution to \eqref{E:wave-ua} provided
\begin{multline}\label{def-weak-sol}
\frac{d^2}{dt^2} \int\limits_{\rN^d \setminus \Omega} u_{\ell}^{a}(t, x) v(x) \, dx + \int\limits_{\rN^d \setminus \Omega} \nabla u_{\ell}^{a}(t, x) \nabla v(x) \, dx \\[6pt]
+ \ell \, \int\limits_{\Gamma} \big(B_{\ell}^\eps \, u_{\ell}^{a} \big)(t, x) \, v(x) \, dx = \int\limits_{\rN^d \setminus \Omega} f(t, x) \, v(x) \, dx \quad \forall \, v \in H^{1} (\rN^{d} \setminus \Omega), 
\end{multline}
for any $t > 0$, and
\begin{equation}\label{def-weak-int}
u_{\ell}^{a}(0, x) =  \partial_t u_{\ell}^{a}(0, x) = 0 \mbox{ in } \rN^{d} \setminus \Omega.
\end{equation}
\end{definition}
\vspace{6 pt}

\begin{remark}
In the definition,  the last term of the LHS of \eqref{def-weak-sol} is 0 if $\ell =0$. The definition in the case $\ell = 0$ is standard. 
\end{remark}

The proof of Theorem~\ref{T:Main} is presented in Section~\ref{sect:Theorem1}. The proof of well posedness of \eqref{E:wave-ua} for $\ell=1$ (non-local structure in time) is based on a nontrivial energy estimate \eqref{est-stability}, which is derived from the causality, see \eqref{E:B}. Following the strategy in \cite{NguyenVogelius1},  we will derive \eqref{k-wave} from estimates in the frequency domain.  For this end, we establish  estimates  for $C(\mf, R)$ in \eqref{known-C} where the dependence on $k$ is well-controlled. This is one of the main parts of the analysis and  presented in the following three propositions which deal with different regimes of frequency.

\begin{proposition}\label{P:LF} Let  $d=2, 3$, $\ell=0, 1$, $0< \eps < 1$,  $0<\mf <\eps^2$, $R_0>0$, and $\ms \in L^2(\rN^d \setminus \overline \Og)$ with $\supp \ms \subset B_{R_0} \setminus  \overline \Og$. Let $v^\eps \in H^1_{\loc}(\mR^d)$ and $v_{\ell}^a \in H^1_{\loc}(\mR^d \setminus \Omega)$ 
be the unique outgoing solutions to \eqref{def-v} and \eqref{def-va} respectively. 
We have, for any $r>0$, 
\begin{equation*}
\|v^{\eps}-v_{\ell}^{a}\|_{H^1(B_r \setminus \Og)} \leq C_r \; \|\ms \|_{L^{2}(\rN^d)},
\end{equation*}
for some  constant $C_r>0$, independent of $\ms$, $\eps$, and $\mf$. 
\end{proposition}

\begin{proposition}\label{P:MF} Let $d=2, 3$,  $\ell=0, 1$, $0< \eps < 1$, $k_0>0$,  $\eps^2 <\mf <k_0$, $R_0>0$, and $\ms \in H^{2\ell+5} (\rN^d \setminus \overline \Og)$ with $\supp \ms \subset B_{R_0} \setminus  \overline \Og$.  Let $v^\eps \in H^1_{\loc}(\mR^d)$ and $v_{\ell}^a \in H^1_{\loc}(\mR^d \setminus \Omega)$ 
be the unique outgoing solutions to \eqref{def-v} and \eqref{def-va} respectively. 
We have, for any $r>0$, \begin{equation}
\|v^{\eps}-v_{\ell}^{a}\|_{H^1(B_r \setminus \Og)} \leq C_r \; \ga^{\ell+1} \; \|\ms \|_{H^{2\ell+5}(\rN^d)},
\end{equation}
for some  constant $C_r>0$, independent of $\ms$, $\eps$, and $\mf$. 
\end{proposition}

\begin{proposition}\label{P:HF}  Let  $d=2, 3$, $\ell=0, 1$, $k_0>0$, $0< \eps < 1$, $\mf\geq k_0$, $R_0>0$, and $\ms \in H^{2\ell+5} (\rN^d \setminus \overline \Og)$ with $\supp \ms \subset B_{R_0} \setminus \overline \Og$. 
Let $v^\eps \in H^1_{\loc}(\mR^d)$ and $v_{\ell}^a \in H^1_{\loc}(\mR^d \setminus \Omega)$ 
be the unique outgoing solutions to \eqref{def-v} and \eqref{def-va} respectively. 
Assume that  $\Og$ is {\bf star-shaped}. Then, for any $K \subset \subset \rN^d \setminus \overline \Og$, we have
\begin{equation}\label{est-P:HF}
\|\nabla(v^\eps-v_{\ell}^{a})\|_{L^2(K)} +\mf \; \|v^\eps-v_{\ell}^{a}\|_{L^2(K)} \leq C_{K} \, \mf^{2\ell+7} \; \ga^{\ell+1}\;  \| \ms \|_{H^{2\ell+5}(\rN^d)},
\end{equation}
for some constant $C_K>0$, independent of  $\ms$, $\eps$, and $\mf$. 
\end{proposition}

The proofs of Propositions~\ref{P:LF} and \ref{P:MF} are given in Section~\ref{sect-Pro12}. They are based on a compactness argument as in \cite{HJNg1} and use results on the Helmholtz equations in the low frequency regime in \cite{Ng-cl-1, Ng-cl-2}. The proof of Proposition~\ref{P:MF} and  \ref{P:HF}  uses  the asymptotic expansion introduced in \cite{HJNg1}. To obtain  explicit dependence on frequency of these estimates, we revise the asymptotic expansion given in \cite{HJNg1} for all range of frequency with a focus on the dependence on frequency. The proof of the stability in Proposition~\ref{P:HF} is given in Section~\ref{sect-Pro3}. It is a heart matter of our paper. The compactness argument used in the proof of  Propositions~\ref{P:LF} and \ref{P:MF} is not appropriate in this regime.  Moreover, the Morawetz's multiplier technique does not work directly in our settings. Due to the structure of the GIBCs, we are only able to obtain an estimate in $L^2(\Gamma)$-norm  of the solution, not the  $H^{1}(\Gamma)$-norm required for Morawetz's technique. To overcome this difficulty,  we use a result  due to H\"ormander in \cite{Hor} (see Lemma~\ref{lem-Hor}). The payoff  for lacking of the control of  the  $H^{1}(\Gamma)$-norm is that we can only obtain estimates in regions away from $\Gamma$, see \eqref{est-P:HF}.

\section{Proof of Theorem \ref{T:Main}} \label{sect:Theorem1}

This section containing two subsections is devoted to the proof of Theorem~\ref{T:Main} assuming Propositions~\ref{P:LF}, \ref{P:MF}, and \ref{P:HF} (their  proofs are given in Section~\ref{sect-Pro12} and \ref{sect-Pro3}). In the first subsection, we establish the well-posedness and the stability for \eqref{E:wave-ua}. We also show that the Fourier transform of the weak solutions satisfies the outgoing conditions for almost every positive frequency. The proof of Theorem~\ref{T:Main} is given in the second subsection.

\subsection{Preliminaries}\label{sect-Theorem1-1}

In this section, we prepare some materials for the proof of Theorem~\ref{T:Main}. We first  prove the well-posedness  and the stability for \eqref{E:wave-ua}. 
 
\begin{lemma}\label{lem-wellposedness}
Let $d=2, 3$, $\ell=0, 1$ and $f \in L^{2}\big([0, \infty), L^{2}(\rN^{d} \setminus \Omega)\big)$ with compact support. There exists a unique weak solution $ v_{\ell}^{a} \in L_{\loc}^\infty\big( [0, \infty), H^{1}(\rN^{d} \setminus \Omega)\big)$ with $\partial_{t} v_{\ell}^{a} \in L_{\loc}^\infty\big([0, \infty), L^{2}(\rN^{d} \setminus \Omega) \big)$ to \eqref{E:wave-ua}.
Moreover, 
\begin{equation}\label{est-stability}
E(t, u_\ell^a)
 \leq C\, t \, \|f\|_{L^{2}([0,t] \times \rN^d)}^{2} \quad \forall \, t  \ge 0. 
\end{equation}
Here,
\begin{equation}\label{energy-1}
E(t, \psi) := \frac{1}{2} \intl_{\rN^d \setminus \Omega} ( |\partial_{t} \psi (t,x)|^{2} + |\nabla \psi(t,x)|^{2} ) \, dx,  
\end{equation}
for $\psi \in L_{\loc}^\infty\big( [0, \infty), H^{1}(\rN^{d} \setminus \Omega)\big)$ with $\partial_{t} \psi \in L_{\loc}^\infty\big([0, \infty), L^{2}(\rN^{d} \setminus \Omega) \big)$. 
\end{lemma}

\noindent{ \bf Proof.} We need only prove the theorem for $\ell=1$, since the case $\ell=0$ is standard.  We first establish the existence of a weak solution  which satisfies \eqref{est-stability}. For this end, we use the Galerkin method. Let $(\varphi_j)_{j=1}^\infty \subset C^\infty_c(\rN^d \setminus \Omega)$ be an orthonormal basis in $H^1(\rN^d \setminus \Omega)$. 
For $m \in \mN$, consider $u_m$ of the form 
\begin{equation}\label{form-uN}
u_m = \sum_{j=1}^{m} d_{m,j}(t) \; \varphi_j(x)
\end{equation}
such that 
\begin{multline}\label{E:GUN}
\frac{d^{2}}{dt^{2}} \intl_{\rN^d \setminus \Omega} u_m(t,x) \; \varphi_j(x) \, dx + \intl_{\rN^d \setminus \Omega} \nabla u_m(t,x)  \; \nabla \varphi_j(x) \, dx \\[6pt] 
+ \int\limits_{\Ga} \big(B_{1}^\eps \, u_m \big)(t,x) \; \varphi_j(x) \, dx = \intl_{\rN^d \setminus \Omega} f(t,x) \; \varphi_j(x) \, dx,\quad \mbox{ for } j=1, \dots, m \end{multline} 
and
\begin{equation} \label{E:intG} 
d_{m,j}(t) = d_{m,j}'(t) = 0 \quad \mbox{ for } j=1, \dots, m.
\end{equation} 
Since $(\varphi_j)_j$ is linearly independent in $H^1(\rN^d \setminus \Omega)$, it is also linearly independent in $L^2(\rN^d \setminus \Omega)$. This implies the $(n \times n)$ matrix $M$ given by $M_{i, j} = \langle \varphi_i, \varphi_j \rangle_{L^2(\rN^d \setminus \Omega)}$ is invertible. 
The existence and uniqueness of $u_m$ then follows; for example, one can use the theory of Volterra equation (see, e.g., \cite[Theorem 2.1.1]{Vol-Burton}).

\medskip
We now derive an estimate for $u_m$. Let us multiply \eqref{E:GUN} by $d_{m,j}'(s)$ and sum it up with respect to $j$.  Integrating the resulting equation over $[0, t]$ with respect to $s$ and using \eqref{E:intG}, we obtain 
\begin{equation}\label{E:Energy1}
E(t, u_{m}) +  \intl_{0}^{t}\intl_{\Ga}\big(B_{1}^\eps \, u_m \big)(s,x) \; \partial_{t} u_{m}(s,x) \, dx \, ds = \intl_{0}^{t}\intl_{\rN^d \setminus \Omega} f(s,x)\; \partial_{t} u_{m}(s,x)\,dx \, ds.
\end{equation}
\medskip
We  claim that, for  $t \ge 0$, 
\begin{equation}\label{claim-pos}
 \intl_{0}^{t}\intl_{\Ga}\big(B_{1}^\eps \, u_m \big)(s,x) \; \partial_{t} u_{m}(s,x) \, dx \, ds \ge 0. 
\end{equation}
Indeed, recall 
\begin{equation*}
\big(B^\eps_1 \, u_m \big) (s,x) = \frac{1} {\sqrt{\pi} \, \eps} \,\Big(\varphi* \pdh_t u_m \Big)(s,x). 
\end{equation*}
where $\varphi$ is given in \eqref{E:conv}.  Here and in what follows $*$ denotes the convolution with respect to time. Set 
\begin{equation*}
U(s,x) = \left\{\begin{array}{cl} \pdh_t u_{m}(s,x) & \mbox{ if } s < t, \\[6pt]
0 & \mbox{ if } s \geq t \mbox{ or } s < 0. 
\end{array}\right. 
\end{equation*}
Then, 
\begin{equation}\label{tototo1}
\intl_{0}^{t} \intl_{\Ga} \big(B^\eps_1 \, u_m \big) (s,x)\; \partial_{t} u_{m}(s,x) \, dx \, ds  = \frac{1}{\eps \sqrt{\pi}} \intl_{\rN} \intl_{\Ga}\Big(\varphi* U \Big)(s,x) \; U(s,x) \, dx \, ds.
\end{equation}
Using Parseval's identity and \eqref{E:conv}, we obtain
\begin{align}\label{tototo2}
\intl_{\rN} \intl_{\Ga} \Big(\varphi* U\Big)(s,x) \; U(s,x) \, dx \, ds
  =&   \intl_{\mR}  \intl_{\Ga}  \widehat{ \Big(\varphi* U\Big)} (\mf,x) \; \overline{\widehat{U}}(\mf,x) \, dx \, d\mf \nonumber \\[6pt] =& 2 \Re \intl_{\mR_+}  \intl_{\Ga} \frac{\sqrt{\pi}}{\ag \, \sqrt{\mf}} \;  |\widehat{U} (\mf,x)|^2 \, dx \, d \mf \geq 0.
\end{align}
A combination of \eqref{tototo1} and \eqref{tototo2} yields  \eqref{claim-pos}. 

\medskip
From \eqref{E:Energy1} and \eqref{claim-pos}, we have
\begin{equation}\label{energy-2}
E(t, u_{m}) \le \intl_{0}^{t}\intl_{\rN^d \setminus \Omega} f(s,x)\; \partial_{t} u_{m}(s,x) \, dx \, ds.
\end{equation}
By the Gronwall inequality, it follows from \eqref{energy-1} and \eqref{energy-2} that 
\begin{equation*}
\int\limits_0^t E(s, u_m) \, ds \le C \left[\int\limits_0^t \Big(\int\limits_0^s \int\limits_{\rN^d \setminus \Omega} |f(\tau,x)|^2 \, dx \, d\tau \Big)^{1/2} \, ds \right]^2, 
\end{equation*}
which implies 
\begin{equation}\label{energy-3}
\int\limits_0^t E(s, u_m) \, ds \le C\, t^2 \, \int\limits_0^t \int\limits_{\rN^d \setminus \Omega} |f(s,x)|^2 \, dx \, ds. 
\end{equation}
Due to \eqref{energy-2},
\begin{equation*}
E(t, u_m)  \le \Big(\int\limits_0^t \int\limits_{\rN^d \setminus \Omega} |\pdh_t u_m(s,x)|^2 \, dx \, ds\Big)^{1/2} \; \Big(\int\limits_0^t \int\limits_{\rN^d \setminus \Omega} |f(s,x)|^2 \, dx \, ds \Big)^{1/2}.
\end{equation*}
This implies 
\begin{equation*}
E(t, u_m)   \le \Big(\int\limits_0^t E(s, u_m) \, ds \Big)^{1/2} \; \Big(\int\limits_0^t \int\limits_{\rN^d \setminus \Omega} | f(s,x)|^2 \, dx \, ds\Big)^{1/2}. 
\end{equation*}
It follows that from  \eqref{energy-3} that
\begin{equation}\label{energy-4}
 E(t, u_m)  \le C t \int\limits_0^t \int\limits_{\rN^d \setminus \Omega} | f(s,x)|^2 \, dx \, ds. 
\end{equation}
Hence, for any fixed $T >0$, there exists a subsequence of $(u_m)$ (which is also denoted by $u_m$ for notational ease) such that $u_m \to u$ weakly-star in $ L^\infty \big([0, T], H^1(\mR^d \setminus \Og) \big)$ and $\partial_t u_m \to \partial_t u $  weakly star in $L^\infty\big([0, T], L^2(\mR^d \setminus \Og) \big)$. It is clear that  $u$ satisfies \eqref{def-weak-sol} for $t \in (0, T)$ and \eqref{def-weak-int}. 

\medskip
\medskip
It remains to show that the limit is unique. It suffices to prove that if $u \in  L^\infty \big([0, T], H^1(\mR^d\setminus \Og ) \big)$ with $\partial_t u  \in L^\infty\big([0, T], L^2(\mR^d \setminus \Og) \big)$ and $u$ satisfies \eqref{def-weak-sol} with $f=0$ and \eqref{def-weak-int} then $u =0$. 

\noindent Set $u^*(t, x) =  \int\limits_0^t u(\tau, x) \, d \tau.$ We claim that $u^*$ satisfies 
\begin{multline}\label{def-weak-sol-*}
\frac{d^2}{dt^2} \int\limits_{\rN^d \setminus \Omega} u^*(t, x) \; v(x) \, dx + \int\limits_{\rN^d \setminus \Omega} \nabla u^*(t, x) \; \nabla v(x) \, dx \\[6pt]
+\int\limits_{\Gamma} \big(B_{1}^\eps \, u^* \big)(t, x) \, v(x) \, dx = 0, \quad \forall \, v \in H^{1} (\rN^{d} \setminus \Omega), 
\end{multline}
for any $t > 0$, and
\begin{equation}\label{def-weak-int-*}
u^*(0, x) =  \partial_t u^*(0, x) = 0 \mbox{ in } \rN^{d} \setminus \Omega.
\end{equation}
The claim (more precisely, equation \eqref{def-weak-sol-*}) can be verified by integrating \eqref{def-weak-sol} with respect to $t$. We basically only need to verify the validity of the boundary term on the LHS. This follows from Lemma~\ref{L:invariant} below. 

From \eqref{def-weak-sol-*}, we derive
\begin{multline}\label{def-weak-sol-**}
 \int\limits_{\rN^d \setminus \Omega} \partial^2_{tt}u^*(t, x) \; v(x) \, dx + \int\limits_{\rN^d \setminus \Omega} \nabla u^*(t, x) \; \nabla v(x) \, dx \\[6pt]
+\int\limits_{\Gamma} \big(B_{1}^\eps \, u^* \big)(t, x) \, v(x) \, dx = 0, \quad \forall \, v \in H^{1} (\rN^{d} \setminus \Omega), 
\end{multline}
Letting $v(x)=u^*_t$ in \eqref{def-weak-sol-**} and integrating in $[0,t]$, we obtain
\begin{equation}\label{ggg}
\intl_0^t \int\limits_{\rN^d \setminus \Omega} u^*_{tt}(t, x) u^*_t(t,x) \, dx + \intl_0^t \int\limits_{\rN^d \setminus \Omega} \nabla u^*(t, x) \nabla u^*_t(t,x) \, dx  + \intl_0^t \int\limits_{\Gamma} \big(B_{1}^\eps \, u^* \big)(t, x) \, u^*_t(t,x) \, dx = 0.
\end{equation}
By the same argument used to obtain \eqref{claim-pos}, we have
$$\intl_0^t \int\limits_{\Gamma} \big(B_{\ell}^\eps \, u^* \big)(t, x) \, u^*_t(t,x) \, dx \geq 0. $$
It follows from \eqref{ggg} that
\begin{equation*}
 E(t, u^*)  \le E(0,u^*) = 0.
\end{equation*}
Therefore, $u^* \equiv 0$ and, hence, $u\equiv 0$. The proof is complete. \proofend

\begin{remark} Similar ideas  are  taken into account  for the proof of the well-posedness of non-local the wave equations in \cite{NguyenVogelius3} in which the Drude-Lorentz model is used to capture the dependence of the material on the frequency.  
\end{remark}

The following lemma which reveals an interesting property of the integral kernel $\varphi$ of $B_1^\eps$ is used in the proof of Lemma~\ref{lem-wellposedness}.
\begin{lemma} \label{L:invariant} Assume that $\psi \in L^\infty_{\loc}[0,\infty)$. Let $$\Psi(t)= \intl_0^t \psi(s) \, ds.$$
Then, for any $t>0$
$$ \intl_0^t (\varphi *\psi)(s) \; ds = (\varphi * \Psi)(t).$$
Here, $\varphi$ is defined in \eqref{E:conv}.
\end{lemma}


\noindent{\bf Proof.} The proof only involves a change of order of integration and an integration by parts. The details are left to the reader.
\proofend

\medskip
Concerning the outgoing condition of $\hat u^a_\ell$ and $\hat u^\eps$, we have the following two results whose proofs are similar to the one  \cite[Theorem A1]{NguyenVogelius2}. The details are left to the reader. 

\begin{lemma}\label{lem-outgoing-2} Let $d=2, 3$, $\ell =0, 1$ and let $\hat{u}_{\ell}^{a}(\mf, x)$ be the Fourier transform of $u_{\ell}^{a}(t,x)$ with respect to $t$. Then, for almost every $k > 0$,  $\hat{u}_{\ell}^{a}(\mf,.) \in H^1_{loc}(\rN^d\setminus \Og)$ is the unique outgoing solution to
\begin{eqnarray*} \left\{\begin{array}{ll} \Delta \hat{u}_{\ell}^{a}(\mf,x) +\mf^2 \; \hat{u}_{\ell}^{a}(\mf,x) =   - \hat{f}(\mf,x), &\hskip 0 pt  \mbox{ in } \rN^d \setminus \Omega,\\[6 pt]
\hat{u}_{\ell}^{a}(\mf,x) + \mD_{\ell}^{\ga}  \; \pn \hat{u}_{\ell}^{a}(\mf,x) =0 & \mbox{ on } \Gamma.\end{array}
\right.
\end{eqnarray*} 
Moreover, \begin{equation*} \mf \; \hat{u}_{\ell}^{a}(\mf,x) \in L^2_{loc}(\rN\times (\rN^d \setminus \Og)).\end{equation*} 
\end{lemma}

\medskip

\begin{lemma}\label{lem-outgoing-3} 
Let $d=2, 3$, $\ell =0, 1$ and let  $\hat{u}^\eps(\mf, x)$ be the Fourier transform of $u^\eps(t,x)$ with respect to $t$. Then, for almost every $k > 0$,  $\hat{u}^\eps(\mf,.) \in H^1_{loc}(\rN^d\setminus \Og)$ is the unique outgoing solution to
\begin{eqnarray*} \Delta \hat{u}^\eps(\mf,x)+ \mf^2 \; \hat{u}^\eps(\mf,x) + i \, \mf \;\sg_\eps \; \hat{u}^\eps(\mf,x) =   - \hat{f}(\mf,x), &\hskip 0 pt  \mbox{ in } \rN^d. 
\end{eqnarray*} 
Moreover,
$$\mf \, \hat{u}^\eps(\mf,x)  \in L^2_{loc}(\rN \times \rN^d).$$
\end{lemma}

\subsection{Proof of Theorem~\ref{T:Main}} \label{sect-Theorem1-2}
The unique existence of $u_{\ell}^a$ and \eqref{est-stability} follow from Lemma \ref{lem-wellposedness}.  It only remains to prove \eqref{k-wave}. To this end, we use Propositions \ref{P:LF}, \ref{P:MF},  and \ref{P:HF} (their proofs are given later in Sections~\ref{sect-Pro12} and \ref{sect-Pro3}). 
We follow the strategy in \cite{NguyenVogelius2}. However, instead of estimating $\|u^{\eps}-u_{\ell}^{a}\|_{L^2(\mR_+, H^1(K))}$ as in \cite{NguyenVogelius2}, we estimate $\|\pdh_{t} (u^{\eps}-u_{\ell}^{a}) \|_{L^2(\mR_+, H^1(K))}.$ This simple idea helps us to avoid the technical issue of integrability in low frequency range in \cite{NguyenVogelius2}, which involves the theory of Gamma-convergence \footnote{This simple idea is also very useful in the context of cloaking in \cite{NguyenVogelius3}.}.


\medskip

Applying  Parseval's identity and using the fact that $u^\eps$ and $u_\ell^a$ are real, we have
\begin{equation} \label{E:Par} \intl_{\rN_+} \|\pdh_{t} (u^{\eps}-u_{\ell}^{a})(t,.) \|^2_{H^1(K)} \;dt = 2 \intl_{\rN_+} \mf^2 \; \|(\hat{u}^{\eps} - \hat{u}_{\ell}^{a})(\mf,.) \|^2_{H^1(K)} \; d\mf.\end{equation}
For notational ease, we assume that the constant $k_0$ in Propositions~\ref{P:MF} and \ref{P:HF} is 1. It is clear that 
\begin{equation}\label{decompose}
 \intl_{\rN_+} \mf^2 \; \|(\hat{u}^{\eps} - \hat{u}_{\ell}^{a})(\mf,.) \|^2_{H^1(K)} \; d\mf  =  \Big( \int\limits_0^{ \eps^2} + \int\limits_{ \eps^2}^1 + \int\limits_{1}^\infty \Big) \mf^2 \; \|(\hat{u}^{\eps} - \hat{u}_{\ell}^{a})(\mf,.) \|^2_{H^1(K)} \; d\mf.
\end{equation}
We next estimate the RHS of \eqref{decompose}. We begin with the first term. Applying Proposition~\ref{P:LF}, we have
\begin{equation*} \intl_{0}^{\eps^2}  \mf^2 \; \|(\hat{u}^{\eps} - \hat{u}_{\ell}^{a})(\mf,.) \|^2_{H^1(K)} \; d\mf \; \leq C\;  \intl_{0}^{\eps^2}  \mf^2 \;  \|\hat{f}(\mf,.)\|^2_{L^2(\rN^d)}  \; d\mf.
\end{equation*}
It follows that 
\begin{equation}\label{est1} \intl_{0}^{\eps^2}  \mf^2 \; \|(\hat{u}^{\eps} - \hat{u}_{\ell}^{a})(\mf,.) \|^2_{H^1(K)} d\mf \;  \leq C\; \eps^{6} \; \sup_{\mf>0}\; \|\hat{f}(\mf,.)\|^2_{L^2(\rN^d)}.
\end{equation}
Since $f(.,x)$ is supported in $[0,T]$,  it follows from the definition of the Fourier transform that 
 $$|\hat{f}(\mf,x)| \leq \frac{1}{\sqrt{2\pi}}\intl_{\rN}|f(t,x)| \;dt \leq C \; \left(\intl_{\rN} |f(t,x)|^2 \, dt\right)^{1/2};$$ 
which implies
\begin{equation}\label{est2}
\|\hat{f}(\mf,.)\|^2_{L^2(\rN^d)} \leq C \intl_{\rN^d} \intl_{\rN} |f(t,x)|^2 \,dt \, dx = C \; \|f\|^2_{L^2(\rN_+ \times \rN^d)}.
\end{equation}
A combination of \eqref{est1} and \eqref{est2} yields 
\begin{equation}\label{E:LF} \intl_{0}^{\eps^2}  \mf^2 \; \|(\hat{u}^{\eps} - \hat{u}_{\ell}^{a})(\mf,.) \|^2_{H^1(K)} \; d\mf  \leq C \, \eps^{6} \, \|f\|^2_{L^2(\rN_+ \times \rN^d)}.\end{equation}
We next estimate the second term of the RHS in \eqref{decompose}. Applying Proposition \ref{P:MF}, we obtain
\begin{equation}\label{est2-1}
 \intl_{\eps^2}^{1}  \mf^2 \; \|(\hat{u}^{\eps} - \hat{u}_{\ell}^{a})(\mf,.) \|^2_{H^1(K)} \; d\mf \le C\; \sup_{\mf>0} \; \|\hat{f}(\mf,.)\|^2_{H^{2\ell+5}(\rN^d)} \intl_{\eps^2}^{1} \mf^2 \;  \ga^{2(\ell+1)}  \; d\mf.
\end{equation}
Similar to \eqref{est2}, we have
\begin{equation}\label{est2-2}
\sup_{\mf > 0} \|\hat{f}(\mf,.)\|_{H^{2\ell+5}(\rN^d)}  \leq C \; \|f\|_{H^{2\ell+5}(\rN_+ \times \rN^d)}.
\end{equation}
Since $\ga^2= \eps^2/\mf$ and $\ell =0, 1$, it follows that 
\begin{equation}\label{est2-3}
 \intl_{\eps^2}^{1} \mf^2 \; \ga^{2(\ell+1)} \; d\mf  = \eps^{2(\ell+1)}  \intl_{\eps^2}^{1} \mf^{1-\ell} \; d\mf   \leq C \; \eps^{2(\ell+1)}.
\end{equation}
A combination of \eqref{est2-1}, \eqref{est2-2}, and \eqref{est2-3} yields 
\ba \label{E:MF} \intl_{\eps^2}^{1}  \mf^2 \; \|(\hat{u}^{\eps} - \hat{u}_{\ell}^{a})(\mf,.) \|^2_{H^1(K)} \; d\mf \leq C \; \eps^{2(\ell+1)} \; \|f\|^2_{H^{2\ell+5}(\rN_+ \times \rN^d)}.\ea
We now estimate the last term of the RHS in \eqref{decompose}.  Applying Proposition \ref{P:HF}, we obtain
\begin{align*} \intl_{1}^\infty \mf^2 \; \|(\hat{u}^{\eps} - \hat{u}_{\ell}^{a})(\mf,.) \|^2_{H^1(K)} \; d\mf  
 \leq &  \intl_{1}^{\infty}\mf^2  \, \ga^{2(\ell+1)} \,  \mf^{2(2\ell+7)}\; \|\hat{f}(\mf,.)\|^2_{H^{2\ell+5}(\rN^d)} \;d\mf \\
 \leq & C \; \eps^{2(\ell+1)} \intl_{1}^{\infty} \, \mf^{3\ell + 15}\; \|\hat{f}(\mf,.)\|^2_{H^{2\ell+5}(\rN^d)} \; d\mf.
 \end{align*}
It follows that 
\begin{equation}\label{E:HF}\intl_{1}^\infty  \mf^2 \; \|(\hat{u}^{\eps} - \hat{u}_{\ell}^{a})(\mf,.)\|^2_{H^1(K)} \; d\mf  \leq C \; \eps^{2(\ell+1)} \|f\|^2_{H^{m_\ell} (\rN_+ \times \rN^d)}, \end{equation}
where \footnote{\label{footnote1} If $\supp f \cap \bar \Omega  = \O$, then $\|\hat{f}(\mf,.)\|^2_{H^{2\ell+5}(\rN^d)}$ can be replaced by $\|\hat{f}(\mf,.)\|^2_{L^2(\rN^d)}$. It follows that $m_\ell$ can be chosen as follows $m_\ell = 8$ if $\ell =0$ and $m_\ell = 9$ if $\ell =1$ and the constant $C$ in \eqref{E:HF} now depends on the distance between $\supp f$ and $\bar \Omega$. } 
 \begin{equation*}
m_\ell = 13 \mbox{ if } \ell =0  \quad \mbox{ and } \quad m_\ell =16 \mbox{ if } \ell =1.
\end{equation*} 
Plugging (\ref{E:LF}), (\ref{E:MF}), and (\ref{E:HF}) into \eqref{decompose}, we obtain:
\begin{equation}\label{final-est}
\intl_{\rN} \mf^2 \; \|(\hat{u}^{\eps} - \hat{u}_{\ell}^{a})(\mf,.) \|^2_{H^1(K)} \; d\mf \leq C \; \eps^{2(\ell+1)} \; \|f\|^2_{H^{m_\ell}(\rN_+ \times \rN^d)}.
\end{equation}
A combination of  (\ref{E:Par}) and \eqref{final-est} yields
\begin{equation}\label{E:pdh} \int\limits_{\rN_+} \|\pdh_{t} (u^{\eps}-u_{\ell}^{a})(t,.)\|^2_{H^1(K)} \leq C \; \eps^{2(\ell+1)} \; \|f\|^2_{H^{m_\ell}(\rN_+ \times \rN^d)}.\end{equation}
Since $u^\eps -  u_{\ell}^a \equiv 0$ at $t=0$, the conclusion follows.  \proofend


\section{Asymptotic expansion for highly conducting obstacle revisited} \label{sect-expansion}

This section is on the asymptotic expansion of $v^\eps$ to \eqref{def-v} with respect to the small parameter $\ga: = \eps/ \sqrt{k}$
and  is essentially based on the work of \cite{HJNg1}. Our goal is to keep track of the frequency dependence there. We recall the notations in \cite{HJNg1} and state estimates which are used in the proof of Propositions~\ref{P:MF} and \ref{P:HF}. Their proofs are given in the appendix. 

\medskip
Define
\begin{equation}\label{def-Omega}
\Og^\delta = \{x \in \Og,~ d(x,\Ga) \leq \delta\}. 
\end{equation}
In what follows, we fix $\delta>0$ small enough such that any $x\in \Og^\delta$ can be written uniquely in the form $x = x_\Gamma + \nu n$, where $(x_\Gamma,\nu) \in \Gamma \times \rN_+$. Here, $\vn$ is the unit normal vector of $\Gamma$ at $x_\Gamma$ pointing {\bf toward} $\Og$.

Let $d=3$,  $\mH$ and $\mG$ be the mean and Gaussian curvatures of $\Ga$ and let ${\cal C}: = \nabla_\Gamma n$ be the curvature tensor on $\Gamma$. Define the tangential operator ${\cal M}$ by the identity
$$
{\cal C} {\cal M} = G I_\Gamma.
$$
One has \cite[(4.4)]{HJNg1}  for $x \in \Omega^\delta$,
\begin{equation}\label{Laplace}
J_{\nu}^3 \Delta = J_\nu \dive_\Gamma (I_\Gamma + \nu {\cal M})^2  \nabla_\Gamma - J_\nu \cdot (I_\Gamma + \nu {\cal M})^2 \nabla_\Gamma + J_\nu^3 \partial_{\nu \nu}^2 +2 J_\nu^2 ({\cal H} + \nu G) \partial_nu,  
\end{equation}
where 
\begin{equation}\label{def-J}
J_\nu: = \det (I + \nu C) = 1 + 2 \nu {\cal H} + \nu^2 G. 
\end{equation}
The following differential operators $\mA_m$ ($1 \le m \le 8$) are defined in \cite{HJNg1} \footnote{The signs in front of $i$ in our formulae are opposite to the ones in \cite{HJNg1}. This is due to the difference between (\ref{def-v}) and  \cite[(2.3)]{HJNg1}, because of different ways to take the Fourier transform in these papers. The reader should keep this fact in mind when comparing our calculations with those in \cite{HJNg1}.}
\begin{align*}
\mA_1 =& 2\mH  \partial_\eta + 6 \eta \mH  (\partial_\eta^2 + i), \\[6pt]
\mA_2 =& \Delta_\Gamma+\mf^2 + 2 \eta (\mG+ 4\mH^2)  \partial_\eta + 3 \eta^2 (\mG+ 4\mH^2) (\partial_\eta^2 +i),\\[6 pt]
\mA_3 =& 2\eta \Big[\mH \Delta_\Gamma + \di_\Gamma (\mM \nabla_\Gamma) -\nabla_\Gamma \mH  \nabla_\Gamma + 3\mf^2 \mH\Big]   + 4 \eta^2 \mH  \Big[(3 \mG + 2\mH^2)  \partial_\eta \Big] \\ \nonumber & + 4 \eta^3 \mH (3\mG + 2 \mH^2) (\partial_\eta^2 +i), \\[6 pt]
\mA_4 =& \eta^2  \Big[\mG \Delta_\Gamma + 4 \mH \di_\Gamma (\mM  \nabla_\Gamma)+ \di_\Gamma (\mM^2  \nabla_\Gamma)\Big] - \eta^2  \Big[\nabla_\Gamma \mG \nabla_\Gamma + 4 \nabla_\Gamma  \mH   (\mM  \nabla_\Gamma) -3\mf^2 (\mG + 4 \mH^2) \Big] \\ \nonumber & +  4 \eta^3 \mG(\mG+4\mH^2) \partial_\eta + 3 \eta^4 \mG(\mG + 4 \mH^2) (\partial_\eta^2 +i), 
\end{align*}
\begin{align*}
\mA_5 =& 2 \eta^3 \Big [\mG  \di_\Gamma (\mM  \nabla_\Gamma) + \mH  \di_\Gamma \Big]  - 2\eta^3 \Big[\nabla_\Gamma \mG (\mM  \nabla_\Gamma) +\nabla_\Gamma \mH  (\mM^2 \nabla_\Gamma) -2\mf^2  \mH  (3\mG + 2\mH^2)\Big] \\ \nonumber & +  10 \eta^4  \mG^2  \mH  \partial_\eta+ 6 \eta^5  \mG^2 \; \mH  (\partial_\eta^2 + i), \\[6pt]
\mA_6 =& \eta^4 \Big[\mG  \di_\Gamma (\mM^2  \nabla_\Gamma) - \nabla_\Gamma  \mG (\mM^2  \nabla_\Gamma ) + 3\mf^2 \mG (\mG + 4 \mH^2) \Big] + 2\eta^5  \mG^3  \partial_\eta+ \eta^6 \mG^3  (\partial_\eta^2 +i),\\[6pt]
\mA_7 =& 6 \eta^5 \mf^2  \mG^3  \mH, \quad \mbox { and } \quad 
\mA_8 = \eta^6 \mf^2 \mG^3. 
\end{align*}
Using \eqref{Laplace} and \eqref{def-J} as in \cite[(5.22)]{HJNg1}, one has 
\begin{equation}\label{IMP}
\Delta + k^2 + \frac{i}{\ga^2} = \frac{1}{J_\nu^3 \ga^2} \Big(-\partial_{\eta \eta}^2  - i  -  \sum_{m=1}^8 \ga^m {\cal A}_m  \Big).  
\end{equation}
Similarly, we also define the above operations in the case $d=2$. In this case,
the triple $(\mH,\mG,\mM)$ is replaced by  $(\frac{\kappa}{2}, 0, 0)$, where $\kappa = \kappa(x)$ is the (signed) curvature of $\Gamma$. 

The following definitions by recurrence of $w_e^\ell$ in $\mR^d \setminus \Omega$ and  $w_i^\ell$ in $\Gamma \times \mR_+$  are given in \cite{HJNg1}. For $\ell = 0$,  define 
\begin{equation}\label{wi0}
w^0_i(x) = 0 \mbox{ in } \Gamma \times \mR_+, 
\end{equation}
and  let $w_e^0 \in H^1_{\loc}(\mR^d \setminus \Omega)$ be the unique outgoing solution to 
\begin{equation}\label{we0} 
\left\{\begin{array}{ll} \Delta w^0_e +\mf^2 \; w^0_e = \ms & \mbox{ in } \rN^d \setminus \overline{\Og}, \\[6pt] 
w_e^0= 0 &  \mbox{ on } \Gamma.
\end{array} \right.
\end{equation}
Let $\ell \ge 1$. Assume that $w_e^j$ and $w_i^j$ are defined for $j \leq \ell-1$. Define $w_i^\ell$ to be the solution to 
\begin{equation} \label{E:Recursive}
\left\{\begin{array}{ll}
(\partial_\eta^2 + i) \; w^\ell_i(x_\Ga,\eta) = - \sum\limits_{m=1}^8  \mA_m \; w^{\ell-m}_i(x_\Ga,\eta) \mbox{ for }  (x_\Ga,\eta) \in \Gamma \times\mR_+, \\[10pt]
\pdh_\eta w_i^\ell(x_\Ga,0) = \pn w_e^{\ell-1}(x_\Ga) \mbox{ and } \lim_{\eta \to \infty} w_i^\ell(x_\Ga,\eta) =0.
\end{array} \right. 
\end{equation}
(here we use the convention $w_i^\ell \equiv 0$ for $\ell<0$) and let  $w_e^\ell \in H^1_{\loc}(\mR^d \setminus \Omega)$ be the unique outgoing solution to 
\begin{equation} \label{E:Recursive-e}
\left\{\begin{array}{cl} \Delta w_e^{\ell} +\mf^2 \; w_e^{\ell} =0  & \mbox{ for } x \in   \rN^d \setminus \overline{\Og}, \\[6pt] w_e^\ell(x)= w_i^\ell(x,0) & \mbox{ for } x \in \Gamma. 
\end{array} \right.
\end{equation}
From \eqref{E:Recursive}, one has, \cite[(4.27) and (4.28)]{HJNg1}, 
\begin{equation}\label{def-wi1}
w_i^1(x_\Gamma, \eta) = - \frac{1}{\alpha }\partial_n w_e^0 \, e^{- \alpha \eta}
\end{equation}
and 
\begin{equation}\label{def-wi2}
w_i^2 (x_\Gamma, \eta) = \left\{  - \frac{1}{\alpha} \, \partial_n w_e^1 + \frac{{\cal H}}{ \alpha^2} \, \partial_n w_e^0 + \frac{\eta {\cal H}}{\alpha}\, \partial_n w_e^0 \right\} e^{-\alpha \eta}. 
\end{equation}
Let $\chi \in C_0^\infty(\rN)$ satisfy \begin{equation} \label{E:Chi} \chi(\eta) = 
\left\{ \begin{array}{cl}  1 & \mbox{ if } |\eta| \leq \delta/2 , \\[6 pt] 
0 & \mbox{ if } |\eta| \geq \delta.\end{array} \right.
\end{equation}
Following \cite{HJNg1}, we set 
\begin{align*} v_{\ell}^{e}(x) =& w_e^0(x) + \ga \, w_e^1 + ...+ \ga^\ell \, w_e^\ell(x),\hskip 140 pt x \in \rN^d \setminus \overline{\Og} \\[6pt] v_\ell^i (x) =& \left[w_i^0 (x_\Gamma, \nu/\ga)+\ga \,  w_i^1 (x_\Gamma,\nu/\ga)+...+  \ga ^{\ell} \, w_i^\ell (x_\Gamma, \nu/\ga)\right] \chi(x),\hskip 20 pt x \in \Og. \end{align*} 
For  $x \in \Og^\delta$,  define
\begin{equation}\label{def-vp}
\varphi_{\ga} (x) = \nu \;\chi(\nu/\ga) \; \partial_n w_e^\ell(x_\Ga) \mbox{ where } x = x_\Gamma + \nu n.
\end{equation}
It is clear that  $\varphi_{\ga} \in C^\infty(\overline \Og)$, 
\begin{equation}\label{pro-vp}
\varphi_{\ga} (x) =0 \mbox{ on } \Gamma \quad  \mbox{ and } \quad \pn \varphi_{\ga} (x) = \pn w_e^\ell \mbox{ on } \Gamma.
\end{equation} 
We define \footnote{Our definition of $\mw_\ell$ is slightly different from the one in \cite{HJNg1} in $\Omega$.  We include the term $-\ga^\ell \varphi_{\ga}(x)$ to make $\pdh_n \mw_\ell$ continuous across $\Gamma$ while maintaining the continuity of $\mw_{\ell}$. This modification is convenient for the use of Morawetz's technique in the high frequency regime later. The scaling for the variable of function $\chi$ in the definition of $\varphi$ reflects the skin effect.} 
 \begin{equation}\label{def-md}
\mw_\ell : = \left\{\begin{array}{ll} v^\eps (x)- v_e^\ell(x)  & \mbox{ in }  \mR^d  \setminus \Omega, \\[6pt]
v^\eps (x)- v_i^\ell(x) -  \ga^\ell \varphi_{\ga}(x) & \mbox{ in } \Omega.
\end{array} \right.
\end{equation}
and 
\begin{equation}\label{E:mq}
\mq_\ell: =  \Delta \mw_\ell +\mf^2\; \mw_\ell + \frac{i}{\ga^2} \; \chi_\Og \; \mw_\ell \quad  \mbox{ in } \mR^d.
\end{equation}
We also set
\begin{equation}\label{def-ed}
\me_\ell := v_{\ell}^{e} - v_{\ell}^{a} \mbox{ in } \mR^d \setminus \Omega \quad \mbox{ and } \quad h_\ell: = \me_\ell + \mD_{\ell}^{\ga} \; \pn \me_\ell \mbox{ on } \Gamma. 
\end{equation}
Then $\me_\ell \in H^1_{\loc}(\mR^d)$ and $\me_\ell$ satisfies the outgoing condition and 
\begin{equation}\label{E:vek} \left\{\begin{array}{ll} \Delta \me_\ell +\mf^2 \; \me_\ell = 0,&\hskip 0 pt \mbox{ in } \rN^d \setminus \overline{\Og}, \\[6 pt] \me_\ell + \mD_{\ell}^{\ga} \; \pn \me_\ell = h_\ell,&\hskip 0 pt \mbox{ on } \Ga,  \end{array} \right.\end{equation}

We now state estimates used later. We begin with two estimates  on $w_e^\ell$. The first one, whose proof is given in Section~\ref{S:ueM}, deals with the low and moderate frequency regimes. 

\begin{lemma}\label{L:ueM} Let $d=2, 3$, $\ell = 0, 1, 2$, $m \geq 0$, $k_0>0$, and $R>0$. There exist two positive constants $C_{R, \ell}$ and $C_{\ell, m}$ independent of $\eps$ and $k$ such that, for  $0< \mf \leq k_0$, 
\begin{equation}\label{E:ueM}
\|w_e^\ell\|_{H^{m+1}(B_R \setminus \Og)}  \leq C_{R, \ell, m} \; \|\ms\|_{H^{2\ell + m}(\rN^d)}
\end{equation}
and, for $m \ge 1$,  
\begin{equation} \label{E:pueM} \|\partial_n w_e^\ell\|_{H^{m-1/2}(\Gamma)} \leq C_{\ell, m}  \; \|\ms\|_{H^{2\ell+m}(\rN^d)}.\end{equation}
\end{lemma}

Here is the  second estimate of $w_e^\ell$ in the high frequency regime whose proof is given in Section~\ref{S:ueH}. 

\begin{lemma}\label{L:ueH} Let $d=2, 3$, $\ell=0, 1, 2$,  $ m \geq 0$, $k_0 > 0$, and $R>0$. Assume that $\Og$ is {\bf star-shaped}. There exist two positive constants $C_{R, \ell}$ and $C_{\ell, m}$ independent of $\eps$ and $k$ such that, for  $ \mf \ge k_0$, 
\begin{equation}\label{E:ueHF}\|w_e^\ell\|_{H^{m+1}(B_R \setminus \Og)} + \mf \; \|w_e^\ell\|_{H^m(B_R \setminus \Og)} \leq C_{R, \ell, m} \, \mf^{2 \ell + m}  \; \|\ms\|_{H^{2 \ell + m}(\rN^d)}
\end{equation}
and, for $m \ge 1$,  
\begin{equation} \label{E:pueHF} \|\partial_n w_e^\ell\|_{H^{m - 1/2}(\Gamma)} \leq C_{\ell, m} \, \mf^{2\ell+m} \|\ms\|_{H^{2\ell+m}(\rN^d)}.
\end{equation}
\end{lemma}
The following two lemmas give us the essential estimates for $\mq_\ell$ and $h_\ell$. Their proofs are given in Sections~\ref{S:wkMH} and \ref{S:ekMH}, respectively. 
\begin{lemma}\label{L:wkMH} We have $\supp \mq_\ell \subset \overline \Og$. Moreover, let $0 < \eps < 1$, $k \ge 0$, $k_0 > 0$,  $\ell = 0, 1, 2$. Assume that $m \ge 1$ and let $s \in H^{\ell + m}(\mR^d)$.  Then, there is a positive constant $C$ independent of $\eps$ and $k$, such that 
\begin{itemize}
\item[i)]  for $\eps^2 < \mf \leq k_0$, we have
 \begin{equation}\label{E:gk} 
\|\mq_\ell\|_{L^2(\rN^d)} \leq C \, \ga^{\ell-1} \; \|\ms\|_{H^{2\ell+3}(\rN^d)}. 
\end{equation}
\item[ii)] assuming in addition that $\Og$ is {\bf star-shaped},  for $\mf \geq k_0$, we have
 \begin{equation}\label{E:gk} 
 \|\mq_\ell\|_{L^2(\rN^d)}  \leq C \, \mf^{2\ell+3} \; \ga^{\ell-1} \; \|\ms\|_{H^{2\ell+3}(\rN^d)},
 \end{equation} 
\end{itemize}
\end{lemma}
\begin{lemma}\label{L:ekMH} Let $0 < \eps < 1$, $k \ge 0$, $k_0 > 0$,  $\ell = 0, 1, 2$. Assume that $m \ge 1$ and let $s \in H^{\ell + m}(\mR^d)$.  Then, there is a positive constant $C_{\ell, m}$ independent of $\eps$ and $k$ such that
\begin{itemize} 
\item[i)] We have, for $\eps^2< \mf \leq k_0$,   
\begin{equation*}
 \|h_\ell\|_{H^{m-1/2} (\Ga)} \leq C_{\ell, m} \, \ga^{\ell+1} \; \|\ms \|_{H^{2\ell+m}(\rN^d)}.
 \end{equation*}
\item[ii)] Assume in addition  that  $\Omega$ is {\bf star-shaped}. We have, for $k \ge k_0$,   
\begin{equation*}
 \|h_\ell\|_{H^{m-1/2}(\Ga)} \leq C_{\ell, m} \, \mf^{2\ell + m} \; \ga^{\ell+1} \; \|\ms \|_{H^{2\ell+m}(\rN^d)}.
 \end{equation*}
\end{itemize}
\end{lemma}


\section{Proofs of Propositions \ref{P:LF} and \ref{P:MF}} \label{sect-Pro12}

This section containing three subsections is devoted to Propositions \ref{P:LF} and \ref{P:MF}. In the first subsection, we present several useful lemmas. The proofs of Propositions \ref{P:LF} and \ref{P:MF} are given  in the last two subsections. 

\subsection{Preliminaries}

In this section, we  present useful lemmas used in the proof of Propositions \ref{P:LF} and \ref{P:MF}. We first recall  the following results established in \cite[Lemmas 2.2]{Ng-cl-2} (see also \cite[Lemma 2.2]{Ng-cl-1}). 

\begin{lemma} \label{lem-pro1-1} 
Let $d=2, 3$,  $0 < k < k_0$, and $\Omega$ be a smooth bounded connected subset of $\rN^d$. Let  $g_\mf \in H^{\frac{1}{2}}(\partial \Og)$ and  $v_\mf \in H^1_{\loc}(\rN^d \setminus \Omega)$ be the unique outgoing solution to 
\begin{equation*}
\left\{\begin{array}{ll}
\Delta v_\mf + \mf^2 \; v_\mf = 0 & \mbox{in } \rN^d \setminus \overline \Omega, \\[6pt]
v_\mf = g_\mf & \mbox{on } \Gamma.
\end{array}\right.
\end{equation*}
Assume that $g_\mf \rightharpoonup g $ weakly in $H^\frac{1}{2}(\partial \Omega)$ as $\mf \to 0$. Then $v_\mf \rightharpoonup v$ weakly in $H^1_{\loc}(\rN^d \setminus  \Omega)$ where $v \in W^1(\rN^d \setminus \bar \Omega)$ is the unique solution to 
\begin{equation*}
\left\{\begin{array}{ll}
\Delta v = 0 & \mbox{in } \rN^d \setminus \overline \Omega, \\[6pt]
v= g & \mbox{on } \Gamma.
\end{array}\right.
\end{equation*}
\end{lemma}

\noindent Here, for an open unbounded subset $U $ of $\rN^d$, the space $W^1(U)$ is defined as follows 
\begin{equation}\label{def-W1}
W^{1}(U) = \left\{ \begin{array}{ll} 
\dsp \Big\{ \psi \in L^1_{loc}(U);~ \frac{\psi(x)}{\ln (2+|x|) \sqrt{1 + |x|^2}}  \in L^2(U) \mbox{ and } \nabla \psi \in L^2(U) \Big\},&\hskip 0 pt \mbox{ if }d=2,~~\\[6 pt]
\dsp \Big\{ \psi \in L^1_{loc}(U);~ \frac{\psi(x)}{\sqrt{1 + |x|^2}}  \in L^2(U) \mbox{ and } \nabla \psi \in L^2(U) \Big\}, &\hskip 0 pt \mbox{ if }d=3.
\end{array}
\right.
\end{equation} 

Using Lemma~\ref{lem-pro1-1}, we can prove

\begin{lemma}\label{lem-stability1}
Let $d=2, 3$, $\ell=0, 1$, $0 < \eps < 1$, $k_0 > 0$, $r_0 > 0$,  $0 < \mf  < k_0$, $q \in L^2(\rN^d \setminus \Omega)$ with $\supp q \subset B_{r_0} \setminus \Omega$, and $g \in H^{1/2}(\Gamma)$. Let $v \in H^1_{\loc}(\rN^d \setminus \Omega)$ be the unique outgoing solution to 
\ba  \label{E:v} \left\{\begin{array}{cl}
\Delta v + \mf^2 \; v = q & \mbox{ in } \rN^d \setminus \Og,\\[6 pt]
v+ \mD_{\ell}^{\ga} \; \pdh_n v = g &  \mbox{ on } \Ga.\end{array}
\right.
\ea
We have 
\begin{equation} \label{E:LM01}
\| v\|_{H^1 (B_r \setminus \Omega)} \le C_r \big( \|q\|_{L^2(\rN^d \setminus \Omega)} + \|g\|_{H^{1/2}(\Gamma)} \big), 
\end{equation}
for some positive constant $C_r$  independent of $g$, $q$, $\eps$, and $\mf$. 
\end{lemma}

\noindent {\bf Proof.}
We only derive the estimate  for small enough $\eps$ and $\mf$. The other case  follows in the same  spirit. 

\medskip
 Set $r_1: = r_0 + 1$. We first prove that, for small enough $\mf$,
\begin{equation}\label{pro1-1}
\|v\|_{L^2(B_{r_1} \setminus \Og)} \leq C \big(\|q\|_{L^2(\rN^d \setminus \Og)} + \| g\|_{H^{1/2}(\Gamma)} \big),
\end{equation}
for some positive constant $C$, independent of $\eps$, $\mf$,  $q$ and $g$, by contradiction. Suppose this is not true. Then there exist $\eps_n \to 0^+$, $\mf_n \to 0^+$, $q_n \in L^2(\rN^d)$ with $\supp q_n \subset B_{r_0} \setminus \Og$, and $g_n \in H^{1/2}(\Gamma)$ such that
\begin{equation}\label{contra-Pro1} 
\|q_n\|_{L^2(\rN^d \setminus \Omega)} + \|g_n \|_{H^{1/2}(\Gamma)} \to 0 \mbox{ and } \|v_n\|_{L^2(B_{r_1} \setminus \Og)}=1. 
\end{equation}
Here, $v_n \in H^1_{\loc}(\rN^d \setminus \Omega) $ is the unique outgoing solution to the problem
\ba \label{E:vN} \left\{\begin{array}{cl}
\Delta v_n + \mf_n^2 \; v_n = q_n & \mbox{ in } \rN^d \setminus \overline{\Og},\\[6 pt]
v_n + \mD_\ell^{\ga_n} \pdh_n v_n =g_n & \mbox{ on } \Ga.\end{array}
\right.
\ea
Using the standard regularity theory of elliptic equations and the representation formula for the equation $\Delta v_n + \mf_n^2  \; v_n = 0$ in $\rN^d \setminus B_{r_0}$, we have
\begin{equation}\label{est-pro1-1}
\| v_n\|_{H^1(B_r \setminus B_{r_0 + 1/2})} \le C_r \mbox{ for all } r > 0. 
\end{equation}
Multiplying the first equation of \eqref{E:vN} by $\overline v_n$ (the conjugate of $v_n$), integrating  over $B_{r_1} \setminus \Omega$, and using the second equation of \eqref{E:vN},  we obtain
\begin{equation}\label{est-pro1-2}
\intl_{B_{r_1} \setminus \Og} \nabla v_n \; \nabla \overline v_n  + \intl_{\Ga} \pn v_n \; \overline{\mD_\ell^{\ga_n}  \pn v_n} = \intl_{\partial B_{r_1}} \pn v_n \; \overline v_n
+  \mf_n^2 \intl_{B_{r_1}\setminus \Og} |v_n|^2  - \intl_{B_{r_1} \setminus \Og} q_n \;  \overline v_n +  \intl_{\Ga} \pn v_n \; \overline{g_n}. 
\end{equation}
Since, by \eqref{notation} and \eqref{E:BFreq},  
\begin{equation}\label{good-sign}
\Re \; \Big[\pn v_n\, \overline{\mD_\ell^{\ga_n} \; \pn v_n}\Big] \ge 0, 
\end{equation}
it follows from \eqref{est-pro1-1} and \eqref{est-pro1-2} that 
\begin{equation}\label{est-pro1-3}
\int\limits_{B_{r_1} \setminus \Omega} |\nabla v_n|^2 \le C,
\end{equation} for a constant $C$ independent of $n$.
From \eqref{contra-Pro1}, \eqref{est-pro1-1}, and \eqref{est-pro1-3},  w.l.o.g. one may assume that $v_n \to v $ weakly in $H^1_{\loc}(\rN^d \setminus \Omega)$ and, by \eqref{E:vN},
\begin{equation}\label{eq1}
\left\{ \begin{array}{cl}
\Delta v  = 0 & \mbox{ in } \rN^d \setminus \overline \Omega, \\[6pt]
v = 0 & \mbox{ on } \Gamma. 
\end{array}\right. 
\end{equation}
Applying Lemma~\ref{lem-pro1-1}, we have 
\begin{equation}\label{v-W1}
 v \in W^1(\rN^d \setminus \bar \Omega). 
 \end{equation}
This implies 
\begin{equation*}
v = 0 \mbox{ in } \rN^d \setminus \Omega. 
\end{equation*}
On the other hand, we derive from \eqref{contra-Pro1} that
 $$\| v\|_{L^2(B_{r_1} \setminus \Omega)} = \lim_{ n \to \infty}\| v_n\|_{L^2(B_{r_1} \setminus \Omega)} =1.$$
We have a contradiction. Thus \eqref{pro1-1} holds.  

\medskip
From \eqref{pro1-1}, as in the proof of \eqref{est-pro1-1} and \eqref{est-pro1-3}, we obtain  \eqref{E:LM01}. The proof is complete. \proofend 

\medskip
We now state the last results in this section dealing with  \eqref{def-v} in the low and moderate frequency  regimes. 

\begin{lemma}\label{lem-stability2} 
Let $d=2, 3$, $0 < \eps < 1$, $k_0 > 0$, $r_0 > 0$,  $0 < \mf  < k_0$, and let  $q \in L^2(\rN^d)$ with $\supp q \subset B_{r_0}$. Let $v \in H^1_{\loc}(\rN^d \setminus \Omega)$ be the unique outgoing solution to 
\begin{equation*}
 \Delta v +\mf^2\; v + i\mf \, \sigma_\eps\, v=  q  \mbox{ in } \rN^d.
 \end{equation*}
We have, for $r>0$, 
\begin{equation}\label{Sta-1}
\|v \|_{H^1(B_r)} \le C_r \big( \|q\|_{L^2(\Omega^c)} + \ga \|q\|_{L^2(\Omega)} \big). 
\end{equation}
where $C_r$ is a positive constant  independent of $\eps$, $\mf$, and $q$. 
\end{lemma}

\medskip
\noindent{\bf Proof.} The proof of Lemma \ref{lem-stability2} is similar to one of Lemma \ref{lem-stability1}. The details are  left to the reader. \proofend

\subsection{Proof of Proposition~\ref{P:LF}} Proposition~\ref{P:LF} is a direct consequence of Lemmas~\ref{lem-stability1} and \ref{lem-stability2}. Indeed, applying Lemma~\ref{lem-stability1} for $v = v_{\ell}^{a}$, $q = \ms$, we obtain, for $\ell =0, 1$, 
\begin{equation}\label{ttt1}
\| v_{\ell}^{a}\|_{H^1 (B_r \setminus \Omega)} \le C \,  \|\ms\|_{L^2(\rN^d)}. 
\end{equation}
and,  applying Lemma~\ref{lem-stability2} for $v = v^\eps$ and $q = \ms$, we have
\begin{equation}\label{ttt2}
\|v^\eps \|_{H^1(B_r)} \le C \; \|\ms\|_{L^2(\rN^d)}.
\end{equation}
A combination of \eqref{ttt1} and \eqref{ttt2} yields the conclusion. \proofend

\subsection{Proof of Proposition \ref{P:MF}}
It is from the definition of $\mw_\ell$ \eqref{def-md} and $\me_\ell$ \eqref{def-ed} that 
\begin{equation} \label{E:vwe} \|v^\eps - v_{\ell}^{a}\|_{H^1(B_r \setminus \Og)} \leq \|\mw_\ell \|_{H^1(B_r \setminus \Og)} + \|\me_\ell \|_{H^1(B_r \setminus \Og)}. \end{equation} 
Applying Lemmas~\ref{lem-stability1}  and \ref{L:ekMH} (with $m=1$), from \eqref{E:vek}, we have for $\ell =0, 1$, 
\begin{equation} \label{E:mel} \|\me_\ell \|_{H^1(B_r \setminus \Og)} \leq  C \; \|h_\ell\|_{H^{1/2}(\Ga)} \leq C \,  \ga^{\ell+1} \; \|\ms\|_{H^{2\ell+1}(\rN^d)}.\end{equation}
Using  Lemmas~\ref{lem-stability2} and \ref{L:wkMH}, from \eqref{E:mq},  we obtain, for $\ell =0,  1, 2$, 
\begin{equation*}  \|\mw_{\ell} \|_{H^1(B_r \setminus \Og)} \leq C \,  \ga \; \|\mq_\ell \|_{L^{2}(\Og)} \leq  C \,  \ga^{\ell} \; \|\ms\|_{H^{2\ell+ 3}(\rN^d)}. \end{equation*}
By Lemma \ref{L:ueM}, it follows that,  for $\ell =0, 1$,  
\begin{align} \label{E:wl+3} \nonumber
\|\mw_{\ell} \|_{H^1(B_r \setminus \Og)} \leq& \ga^{\ell+1} \; \|w_e^{\ell+1}\|_{H^1(B_r \setminus \Og)} + C \,  \ga^{\ell + 1} \; \|\ms\|_{H^{2\ell+ 5}(\rN^d)} \\[6 pt]  \leq& C \,  \ga^{\ell + 1} \; \|\ms\|_{H^{2\ell+ 5}(\rN^d)} . 
\end{align}
A combination of \eqref{E:vwe}, \eqref{E:mel}, and \eqref{E:wl+3} yields the conclution. \proofend

\section{Proof of Proposition~\ref{P:HF}}\label{sect-Pro3}

This section is devoted to the proof of Proposition~\ref{P:HF}. In order to obtain the desired estimate for $u_{a}^{\ell}-u^\eps$, we will derive separate estimates for the functions $\mw_{\ell}$ and $\me_{\ell}$, introduced in Section \ref{sect-expansion}. This goal is achieved by Corollaries \ref{cor-high} and \ref{cor-high2} below. Our presentation is divided into two subsections. In the first one, we present some useful lemmas. The proof of Proposition~\ref{P:HF} is given in the second subsection. 

\subsection{Preliminaries}

In this section, we present useful lemmas used in the proof of Proposition~\ref{P:HF}. 
We start this section with the following lemma. 

\begin{lemma}\label{L:NV} Let $d=2, 3$, $D$ be a {\bf star-shaped}  domain of $\mR^d$,  and ${r_*}>0$ such that $D \subset B_{r_*}$. Define
\begin{equation*}
P(r) =
\left\{
\begin{array}{cl}
\dsp \frac{2 r_*}{d-1} & \mbox{ if }~ r > {r_*}~,\\[6pt]
\dsp \frac{2r}{d-1} & \mbox{ if }~ 0< r < {r_*}~,
\end{array}\right.
\quad \mbox{and} \quad Q(r) =
\left\{
\begin{array}{cl}
\dsp \frac{r_*}{r}  & \mbox{ if }~ r > {r_*}~,\\[6pt]
\dsp 1 & \mbox{ if }~ 0 < r < {r_*}~,
\end{array}\right.
\end{equation*}
and let  $v \in H^1_{\loc}(\rN^d)$ be such that $\Delta v +\mf^2 v \in L^2_{\loc}(\rN^d)$. For any $R> {r_*}$ and $\mf > 0$, we have
\begin{multline}\label{est-rem}
 \Re \intl_{B_R \setminus \Omega} \big(\Delta v +\mf^2\; v \big)  \big[P(r) \, \overline v_r + Q(r) \; \overline v \big]
 \le \\ - \frac{1}{d-1} \intl_{B_{r_*} \setminus \Omega} \Big( |\nabla v|^2 +\mf^2\; |v|^2 \Big) +\frac{r_* \, (3-d)}{2}  \intl_{B_R \setminus B_{r_*}} \frac{u^2}{r^3}  +  F_0(v) -F(R,v).
\end{multline}
Here, $n$ denotes the {\bf inward} unit normal vector of $\pdh D$, 
\begin{equation*}
F_0(v) = \Re \; \Big( \intl_{\pdh D}  \frac{2}{d-1} \pn v \; (x \cdot \nabla \overline v)   -  \frac{1}{d-1}  (x \cdot \vn) \; |\nabla v|^2+  \pn v \; \overline v +  \frac{\mf^2}{d-1}   (x \cdot \vn) \; |v|^2  \Big),
\end{equation*}
and 
\begin{equation*} F(r,v) = \Re \; \Big(-\intl_{\partial B_r}  \frac{ {r_*}}{d-1} |v_r|^2  - \frac{{r_*}}{d-1}  |\nabla_{\partial B_r} v|^2 + \frac{{r_*}}{2r^2}  |v|^2 + \frac{{r_*}}{r} \, v_r \, \bar v  + \frac{\mf^2 {r_*}}{d-1} |v|^2 \Big).\end{equation*} 
\end{lemma}

\begin{remark} This lemma 
 has been stated and proved in  \cite[Lemma 1]{NguyenVogelius1} for the spherical domains. The proof presented here follows heavily from the one of \cite[Lemma 1]{NguyenVogelius1}. In the proof, we also use the  ``Rellich" identity  \eqref{ident} which has root  from \cite{Rellich, PayneWeinberge, MorawetzLudwig}. 
Estimate \eqref{est-rem} is in the spirit of Morawetz-Ludwig \cite{MorawetzLudwig}. The choice of the weight functions $P(r), Q(r)$ appeared in the work of Perthame and Vega \cite{PerthameVega}.  
\end{remark}

\noindent {\bf Proof.} We will prove the lemma for $v \in C^\infty(\rN^d)$. The general case follows by a standard regularizing argument. 
We have
\begin{equation*}\label{E:brog} \Re \intl_{B_R \setminus D}\big(\Delta v +\mf^2 v \big)\;  \big[P(r) \, \overline v_r + Q(r) \, \overline v \big]= A_1 + A_2,  \end{equation*}
where
\begin{equation*} A_1 =  \Re \intl_{B_{r_*} \setminus D} \big(\Delta v +\mf^2 v \big)\; \big[P(r) \; \overline  v_r + Q(r) \, \overline v \big]   \end{equation*}
and 
\begin{equation*} A_2 =  \Re \intl_{B_R \setminus B_{r_*}}  \big(\Delta v +\mf^2 v \big) \big[P(r) \, \overline v_r + Q(r) \, \overline v \big].  \end{equation*}

\noindent{\underline{Calculate $A_1$}:}  Since $P(r) = \frac{2r}{d-1}$ and $Q(r)=1$ for $0 < r <  {r_*}$,
\begin{equation*} 
A_1 = \Re \intl_{B_{r_*} \setminus D}\big(\Delta v +\mf^2 v \big)  \Big[\frac{2r}{d-1} \overline v_r + \overline v \Big] =  \intl_{B_{r_*} \setminus D}\Re \Big[\big( \Delta v +\mf^2 v \big)  \Big( \frac{2}{d-1} \, x \cdot \nabla \overline v + \overline v \Big)\Big].  
\end{equation*}
We have \footnote{This is the ``Rellich" identity  \eqref{ident} which has root  from \cite{Rellich, PayneWeinberge, MorawetzLudwig}.}
\begin{multline} \label{ident} \Re \; \Big[\left(\Delta v +\mf^2 v\right) \Big(\frac{2}{d-1}x \cdot \nabla \overline v+ \overline v \Big) \Big]  =  - \frac{1}{d-1} \big(|\nabla v|^2 +\mf^2 |v|^2\big) \\ +  \Re \; \nabla \; \Big[\frac{2}{d-1} \nabla v \; (x \cdot \nabla \overline v) -  \frac{1}{d-1}  x \; |\nabla v|^2  + \nabla v \; \overline v + \frac{\mf^2}{d-1}x\; |v|^2 \Big].  \end{multline}
Integrating over the domain $B_{r_*} \setminus D$, we obtain:
\begin{align*}A_1 =&- \frac{1}{d-1} \intl_{B_{r_*} \setminus D} \big(|\nabla v|^2 +\mf^2|v|^2\big) +  \Re \; \intl_{\pdh B_{r_*}}  \frac{2 r_*}{d-1}  \; |v_r|^2  -  \frac{{r_*}}{d-1} \; \big|\nabla v \big|^2 + v_r \; \overline v + \frac{\mf^2 \; {r_*}}{d-1} \; |v|^2  \\ & +  \Re \; \intl_{\pdh D}  \frac{2}{d-1} \pn v \; (x \cdot \nabla \overline v)  -  \frac{1}{d-1} \left(x \cdot \vn\right) \, |\nabla v|^2 + \pn v\; \overline v + \frac{\mf^2}{d-1} \left(x \cdot \vn\right) |v|^2.\end{align*}
It follows that 
\begin{multline*} A_1 =  - \frac{1}{d-1} \intl_{B_{r_*} \setminus D} \big(|\nabla v|^2 +\mf^2|v|^2\big) + F_0(v) + \Re \; \intl_{\pdh B_{r_*}}  \frac{2 r_*}{d-1} |v_r|^2  -  \frac{{r_*}}{d-1} \; |\nabla v|^2 +  v_r \; \overline v + \frac{\mf^2 \, {r_*}}{d-1}  \; |v|^2.  \end{multline*}
Since $|\nabla v|^2  = |v_r|^2 + |\nabla_{\pdh B_{r_*}}v|^2$, we have
\begin{multline*}  \Re \; \intl_{\pdh B_{r_*}}  \frac{2 r_*}{d-1} |v_r|^2  -  \frac{{r_*}}{d-1} \; |\nabla v|^2 + v_r \; \overline v +  \frac{\mf^2 \, {r_*}}{d-1}  \; |v|^2 \\[6pt]
=  \Re \; \intl_{\pdh B_{r_*}}  \frac{{r_*}}{d-1} \; |v_r|^2  -  \frac{{r_*}}{d-1} \; |\nabla_{\pdh B_{r_*}} v|^2 +  v_r \; \overline v +   \frac{{r_*}\mf^2}{d-1} \; |v|^2\le  - F(r_*,v).
\end{multline*}
It follows that
\begin{equation}\label{tt1}  A_1 \leq  - \frac{1}{d-1} \intl_{B_{r_*} \setminus D} \big(|\nabla v|^2 +\mf^2|v|^2\big) + F_0(v) - F(r_*,v).  \end{equation}

\noindent{\underline{Estimate $A_2$}:} Applying  \cite[Lemma 2]{NguyenVogelius1}, we have 
\begin{equation}\label{tt2}
A_2 \leq  \frac{r_* \, (3-d)}{2}  \intl_{B_R \setminus B_{r_*}} \frac{u^2}{r^3}  + F({r_*}, v) - F(R,v).
\end{equation}
The conclusion  now follows from  \eqref{tt1} and \eqref{tt2}. \proofend

\medskip
The following lemma, in spirit of Morawetz-Ludwig \cite{MorawetzLudwig} (see also \cite{PerthameVega}), is important for our analysis. 

\begin{lemma} \label{L:ExtHF} Let  $d=2, 3$, $k_0 > 0$, $r_0 >0$, $q\in L^2(\rN^d)$ with $\supp q \subset B_{r_0} \setminus \Og$, and  $g \in H^{1}(\Ga)$. Let $k \ge k_0$ and $v \in H^1_{\loc}(\rN^d \setminus \Og)$ be the unique outgoing solution to
\begin{equation}\label{E:ext}
\left\{\begin{array}{ll}
\Delta v +\mf^2 \; v = q & \mbox{ in } \rN^d \setminus \overline \Og, \\[6pt]
v = g &  \mbox{ on } \Ga.
\end{array}\right.
\end{equation}
Assume that $\Og$ is {\bf star-shaped}. Given $r_*>0$, there exists a positive constant $C=C({r_*}, r_0 ,k_0, \Og)$ independent of $k$ such that 
\begin{equation}\label{Mor-Lud}
 \|\nabla v\|_{L^2(B_{r_*} \setminus \Og)} +\mf \; \|v\|_{L^2(B_{r_*} \setminus \Og)} + \|\pn v\|_{L^2(\Gamma)}  \leq  C \left(\|q\|_{L^2(\rN^d  \setminus \Og)} +  \|\nabla_\Gamma g\|_{L^2(\Gamma)}+\mf \; \|g\|_{L^2(\Gamma)}\right).
\end{equation} 
\end{lemma}

\noindent{\bf  Proof of Lemma~\ref{L:ExtHF}.} The idea is to apply Lemma \ref{L:NV} for $D= \Og$. We have
\begin{multline*} F_0(v) = \Re \int\limits_{\pdh \Og} \Big( \frac{1}{d-1} (x \cdot n) |\pdh_n v|^2 + \frac{2}{d-1} \pdh_n  v \, \nabla_{\Ga} \bar v - \frac{1}{d-1} (x \cdot n) |\nabla_{\Ga} v|^2  \\+ \pdh_n v \, \bar v + \frac{\mf^2}{d-1}(x \cdot n) \,|v|^2 \Big), 
\end{multline*}
where $n$ denotes the {\bf inward} unit normal vector of $\pdh \Omega$.  Since $\Og$ is star-shaped, it follows that 
\begin{equation}\label{F0} F_0(v)  \leq - \frac{1}{C} \; \|\pdh_n v\|^2 + C \; \left(\|\nabla_\Ga v\|_{L^2(\Ga)} + \mf^2 \|v\|_{L^2(\Ga)} \right).
\end{equation}
Applying Lemma~\ref{L:NV} and using \eqref{F0}, we obtain 
\begin{multline*}
 \|\nabla v\|^2_{L^2(B_{r_*} \setminus \Og)} +\mf^2 \; \|v\|^2_{L^2(B_{r_*} \setminus \Og)} + \|\pn v\|^2_{L^2(\Gamma)} \\ \leq  C \Big( \frac{r_* \, (3-d)}{2}  \intl_{B_R \setminus B_{r_*}} \frac{u^2}{r^3}    + \|q\|^2_{L^2(\rN^d \setminus \Og)} +  \|\nabla_\Gamma g\|^2_{L^2(\Gamma)}+\mf \; \|g\|^2_{L^2(\Gamma)}\Big),
\end{multline*} 
This  implies the conclusion in the case $d=3$. For $d=2$, it remains to absorb the first term in the RHS into the LHS. Without loss of generality, we may assume that $r_*$ is big enough. The absorption then can be done as in \cite[p. 11-12]{NguyenVogelius1}.  The details are left to the reader. \proofend
\medskip

When the control is only available on $L^2(\Gamma)$ (not $H^1(\Gamma)$), one has the following result by H\"ormander \cite[Theorem 3.1]{Hor} (see also \cite[page 65]{Hor}). 

\begin{lemma}[H\"ormander] \label{lem-Hor} Let $D$ be a bounded smooth domain of $\mR^d$ ($d \ge 2$) and $g \in H^{1/2}(\partial D)$. Assume $v \in H^1(D)$ is the unique solution to the system
\begin{equation}
\left\{\begin{array}{ll}
\Delta v = 0 & \mbox{ in } D, \\[6pt]
v = g & \mbox{ on } \partial D. 
\end{array} \right.
\end{equation}
Then 
\begin{equation}\label{est-Hor}
\| v\|_{L^2(D)} \le C \| g\|_{L^2(\partial D)}, 
\end{equation}
for some positive constant $C$ independent of $g$ \footnote{In \eqref{est-Hor},  $\| g\|_{L^2(\partial D)}$ is used not $\| g\|_{H^{1/2}(\partial D)}$}. 
\end{lemma}

Here is a result related to equation \eqref{E:vek} of $\me_\ell$:

\begin{lemma}\label{lem-pro3-1} Let $d=2, 3$, $\ell=0,1$, $0<\eps<1$, $k_0> 0$, $\mf \geq k_0$,   $h \in L^{2}(\Gamma)$, and let $v \in H^1_{\loc}(\rN^d \setminus \Omega)$ be the unique outgoing solution to 
\begin{equation}\label{est-lem-pro2-1}
\left\{\begin{array}{ll}  \Delta v +\mf^2\; v =  0 & \mbox{ in } \rN^d \setminus \overline{\Og}, \\[6pt] v + \mD_{\ell}^{\ga}\; \pn v= h &\mbox{ on } \Gamma. \end{array} \right.
\end{equation}
We have
\begin{equation*}
\|v \|_{L^2(\Gamma)} \le C \; \|h\|_{L^2(\Gamma)}, 
\end{equation*}
where $C$ is a positive constant, independent of $\eps$, $\mf$, $q$, and $h$. 
\end{lemma}

\noindent{\bf Proof.} We only  consider  the case $\ell=1$ since the lemma is trivial for $\ell=0$. Let $\ell =1$.  Multiplying  the first equation of \eqref{est-lem-pro2-1} by $\bar v$, integrating in $B_r \setminus \Og$, and using the boundary condition, we obtain 
\begin{equation}\label{est1-lem-2}
\int\limits_{B_r \setminus \overline{\Og}} -|\nabla v|^2 + \mf^2 \; |v|^2  \;+  \intl_{\pdh B_r} \pdh_r v\; \overline{v} +  \int\limits_{\Gamma} \pn v \;  \overline{h} - \int\limits_{\Gamma} \pn v \;  \overline{\mD_1^{\ga} \, \pn v} =0. 
\end{equation}
Recall that 
\begin{equation} \label{est11-lem-2} 
\mD_1^{\ga} = \frac{\ga}{\alpha} = \frac{\sqrt{2}}{2} \ga + i \frac{\sqrt{2}}{2} \ga,  
\end{equation}
which implies 
\begin{equation}\label{est2-lem-2} 
- \Im \; \Big[ \pn v \; \overline{ \mD_\ell^{\ga} \,  \pn v} \Big] = \frac{\sqrt{2}}{2} \ga \; |\pn \,v|^2. 
\end{equation}
Since $v$ satisfies the outgoing condition \eqref{OGC}, i.e.,  
\begin{equation*}
\partial_r v -  ik v= o(r^{-(d-1)/2}), 
\end{equation*} 
it follows that 
\begin{equation}\label{est3-lem-2}
\liminf_{ r \to \infty }\Im \Big( \intl_{\pdh B_r} \pdh_r v\; \overline{v} \Big) =  \liminf_{r \to \infty} \intl_{\partial B_r} k|v|^2 \ge 0. 
\end{equation}
Considering the imaginary part of \eqref{est1-lem-2} and letting $r \to \infty$, we derive from  \eqref{est2-lem-2} and \eqref{est3-lem-2} that 
\begin{equation*}
\int\limits_{\Gamma} \ga \; |\pn v|^2 \le C \int\limits_{\Gamma} |h| |\pn v|.
\end{equation*}
This implies 
\begin{equation}\label{est4-lem-2}
\ga \; \|\pn v \|_{L^2(\Gamma)} \le \| h\|_{L^2(\Gamma)}. 
\end{equation}
Since
$v =  h - \mD_{\ell}^{\ga} \; \pn v$, the conclusion follows from \eqref{est11-lem-2} and \eqref{est4-lem-2}. \proofend

\medskip

Here is an important consequence of Lemmas \ref{L:ExtHF}, \ref{lem-Hor}, and \ref{lem-pro3-1}, which will be applied to obtain the estimate for $\me_{\ell}$:

\begin{corollary}\label{cor-high}
Let $\ell=0, 1$, $k_0>0$,  $0 < \eps < 1$, $\mf \geq k_0$,  $h \in H^{1/2}(\Gamma)$, and let $v \in H^1_{\loc}(\rN^d \setminus \Omega)$ be the unique outgoing solution to 
\begin{equation}\label{E:cor-high}
\left\{\begin{array}{cl}  \Delta v +\mf^2\; v = 0 & \mbox{ in } \rN^d \setminus \overline{\Og}, \\[6pt] 
v + \mD_{\ell}^{\ga}\; \pn v = h &\mbox{ on } \Gamma. \end{array} \right.
\end{equation}
Assume that $\Og$ is {\bf star-shaped}.  Then, for all $K \subset \subset \rN^d \setminus \overline \Omega$, 
\begin{equation*}
\|\nabla v\|_{L^2(K)} +\mf \; \|v\|_{L^2(K)} \le C \; \mf^2 \; \|h\|_{L^2(\Ga)},
\end{equation*}
for some positive constant $C = C(K)$ independent of $\eps$, $\mf$, and $h$. 
\end{corollary}

\noindent{\bf Proof.} Let $r>0$ such that $K \subset B_r \setminus \Og$. Set $r_1 = r + 1$ and let $\phi \in H^1(B_{r_1} \setminus \Og)$ be the solution to
\begin{equation}\label{E:phi}
\left\{\begin{array}{l}
\Delta \phi = 0  \mbox{ in } B_{r_1} \setminus \bar \Omega, \\[6pt] 
\phi = v   \mbox{ on  } \Gamma \quad \mbox{ and } \quad  
\phi = 0  \mbox{ on } \partial B_{r_1}. 
\end{array} \right.
\end{equation}
Applying  Lemma~\ref{lem-pro3-1} we have
\begin{equation*}
\| v\|_{L^2(\Gamma)} \le C \| h\|_{L^2(\Gamma)}. 
\end{equation*}
It follows from Lemma~\ref{lem-Hor} that
\begin{equation}\label{t2}
\| \phi \|_{L^2(B_{r_1} \setminus \Omega)} \le C \| h\|_{L^2(\Gamma)}. 
\end{equation}
Fix $\chi \in C^\infty (\rN^d)$ such that $\chi = 1$ in $B_{r}$ and $\supp \chi \subset B_{r + 1/2}$. Set 
\begin{equation}\label{def-V}
V = v - \chi \phi \mbox{ in } \rN^d \setminus \Omega. 
\end{equation}
It is clear that $V \in H^1_{\loc}(\rN^d \setminus \Omega)$ is the unique outgoing solution to the problem
\begin{equation}\label{eq-Vk}\left\{
\begin{array}{ll}
\Delta V + \mf^2 \; V = - \Delta (\chi \phi) - \mf^2 \; \chi \phi & \mbox{ in } \rN^d \setminus \overline \Omega, \\[6pt]
V = 0 & \mbox{ on }  \Gamma.  
\end{array} \right.
\end{equation}
Since $\Delta \phi = 0$ in $B_{r_1} \setminus \Omega$, $\chi = 1$ in $B_{r}$, and $\chi = 0$ in  $\rN^d \setminus B_{r  + 1/2}$, 
\begin{equation*}
 \|  \Delta(\chi \phi) \|_{L^2(\rN^d \setminus \Omega)} + \mf^2 \; \| \chi \phi \|_{L^2(\rN^d \setminus \Omega) } \le C \; \big(\|\phi\|_{H^2(B_{r+1/2} \setminus B_r)} + \mf^2 \; \|\phi\|_{L^2(B_{r_1} \setminus \Og)} \big).
\end{equation*}
Using the standard regularity of elliptic equations, we derive from \eqref{E:phi} that 
\begin{equation}\label{est-Vk}
\|  \Delta(\chi \phi) \|_{L^2(\rN^d \setminus \Omega)} + \mf^2 \; \| \chi \phi \|_{L^2(\rN^d \setminus \Omega) } \le C \; \mf^2 \;  \|\phi\|_{L^2(B_{r_1} \setminus \Og)}.
\end{equation}
Applying Lemma~\ref{L:ExtHF}  for  \eqref{eq-Vk} and using \eqref{est-Vk}, we arrive to:
\begin{equation}\label{t3} 
\| \nabla V\|_{L^2(B_r \setminus \Omega)} + \mf \; \| V\|_{L^2(B_r \setminus \Omega)} \le C\; \mf^2 \; \| h\|_{L^2(\Gamma)}. 
\end{equation}
The conclusion now follows from \eqref{E:phi}, \eqref{def-V}, \eqref{t3}, and the standard regularity theory for elliptic equations.
\proofend

\medskip

\noindent The following lemma, which is a variant of  \cite[Proposition 1]{NguyenVogelius2}, plays an important role in analyzing $\mw_\ell$.  

\begin{lemma}\label{L:nhHF} Let $d=2, 3$, $k_0> 0$, $r_0>0$, $0 < \eps < 1$, $\mf \geq k_0$, and $\mq \in L^2(\rN^d)$ with $\supp \mq \subset \Og$.  Let $v \in H^1_{\loc}(\rN^d)$ be the unique outgoing solution to 
\begin{equation} \label{E:nh} \Delta v+\mf^2\; v + i \,\mf \, \sigma_\eps \; v = \mq \mbox{ in } \rN^d.
\end{equation} 
We have, for $r_*>0$,
\begin{equation}\label{H:p1}
\|\nabla v \|_{L^2(B_{r_*})} + \mf  \|v\|_{L^2(B_{r_*})}  \leq C_{r_*} \; \| q\|_{L^2(\Omega)} 
\end{equation}
and 
\begin{equation}\label{H:p2}
\|v \|_{L^2(\Omega)} \le C \, \eps^2\, \|q \|_{L^2(\Omega)} 
\end{equation}
where $C_{r_*}$ and $C$ are positive constants  independent of $\mf$, $\eps$, and $q$. As a consequence,
\begin{equation}\label{H:p3}
\|v \|_{L^2(\Gamma)} \le C \, \eps \, \|q \|_{L^2(\Omega)}.
\end{equation}
\end{lemma}

\noindent {\bf Proof.} We follow the strategy used in the proof  of \cite[Proposition 1]{NguyenVogelius2}. Multiplying equation \eqref{E:nh} by $\overline{v}$ and integrating over $B_R$, we have
\begin{equation*}
-\int\limits_{B_R} |\nabla v|^2 +\mf^2 \int\limits_{B_R} |v|^2 + \int\limits_{\partial B_R} \partial_r v \; \bar v + i k \int\limits_{\Omega} \sigma_\eps \;  |v|^2 = \int\limits_{\Omega} q \; \bar v. 
\end{equation*}
Letting $R \to \infty$, using the outgoing condition, and considering the imaginary part, we obtain 
\begin{equation}\label{important}
\mf\; \limsup_{R \to \infty} \int\limits_{\partial B_R} |v|^2 + \frac{1}{\ga^2} \int\limits_{\Og} |v|^2 \leq \|v\|_{L^2(\Og)} \; \|\mq\|_{L^2(\Og)}.
\end{equation}
This implies
\begin{equation} \label{E:Inf}
\int\limits_{\Og} |v|^2  \leq \ga^4 \int\limits_{\Og} |q|^2 \quad \mbox{ and } \quad\mf\; \limsup_{R \to \infty} \int\limits_{\partial B_R} |v|^2 \leq \ga^2 \int\limits_{\Og} |q|^2.\end{equation}
Hence \eqref{H:p2} is proved. Let $\gamma_0>0$ such that $B_{2 \gamma_0} \subset \Og$. It only remains to prove \eqref{H:p1}. 

\medskip

Multiplying  (\ref{E:nh}) by $\phi^2 \bar v$, with $\phi \in C^\infty_{\mc}(B_{2 \gamma_0})$ and $\phi = 1$ in $B_{\gamma_0}$, and integrating over $\Og$. We obtain, by Caccioppoli's inequality, 
\begin{equation*}
\int\limits_{B_{\gamma_0}} |\nabla v|^2   \le C \;\Big( \mf^2  \int\limits_{\Omega} |v|^2 +  \int\limits_{\Omega} |q|^2 \Big).  
\end{equation*} 
In this proof, $C$ denotes a positive constant independent of $\eps$, $k$, and $q$. It follows from \eqref{E:Inf} that
\begin{equation*}
\int\limits_{B_{\gamma_0}} |\nabla v|^2 +\mf^2 \int\limits_{B_{\gamma_0}} |v|^2 \le C (\ga^4\mf^2+1)  \int\limits_{\Omega} |q|^2 = C (\eps^4+1)  \int\limits_{\Omega} |q|^2 \leq C  \int\limits_{\Omega} |q|^2.
\end{equation*}
Hence, there exists $\gamma_0/2<\tau \leq \gamma_0$ such that 
\begin{equation}\label{HF-0}
\int\limits_{\partial B_{\tau}} |\nabla v|^2 +\mf^2 \int\limits_{\partial B_{\tau}} |v|^2 \le C  \; \int\limits_{\Omega} |q|^2.
\end{equation}
 Applying Lemma \ref{L:NV} with  $D = B_\tau$, we obtain, for any $R> r_* > R_0$, 
\begin{multline}  \label{HF-1} 
 \Re \intl_{B_R \setminus B_\tau} \big[  r \overline v_r +  \overline v \big] \; \big[\Delta v +\mf^2 \; v \Big] 
 \le    - \frac{1}{d-1} \intl_{B_{r_*} \setminus B_\tau} \Big( |\nabla v|^2 +\mf^2 \; |v|^2 \Big)  \\[6pt]
  + \frac{r (3-d)}{2} \int\limits_{B_R \setminus B_{r_*}} \frac{u^2}{r^3}  +  F_0(v) - F(R, v). 
 \end{multline}
Here $F_0(v)$ and $F(R, v)$ are defined in Lemma~\ref{L:NV}.
Using (\ref{E:nh}), we derive from \eqref{HF-1} that
\begin{multline}\label{HF-2} 
 \frac{1}{d-1} \intl_{B_{r_*} \setminus B_\tau} \Big( |\nabla v|^2 +\mf^2 \; |v|^2 \Big) 
 \leq  -\Re \intl_{\Og \setminus B_\tau} \left[r \overline v_r +  \overline v \right] \left[-\frac{i}{\ga^2} v +  \mq\right]  \\[6pt]
  +  \frac{r_* (3-d)}{2} \int\limits_{B_R \setminus B_{r_*}} \frac{u^2}{r^3}
  + F_0(v) - F(R,v).
\end{multline}
We have
\begin{equation*}\label{HF-3}
\left|\Re \intl_{\Og \setminus B_\tau} \left[r  \overline v_r + \overline v \right] \; \left[-\frac{i}{\ga^2} v +  \mq \right] \right| \le  \frac{C}{\ga^2}  \intl_{\Og \setminus B_\tau} |v_r| \; |v|  + C \intl_{\Og \setminus B_\tau} (|v_r| + |v| ) |\mq|.
\end{equation*}
We derive from Young's inequality that 
\begin{equation}\label{HF-5}
\left|\Re \intl_{\Og \setminus B_\tau} \left[r  \overline v_r + \overline v \right] \; \left[-\frac{i}{\ga^2} v +  \mq \right] \right|  \le \frac
{1}{2(d-1)} \intl_{\Og \setminus B_\tau} |v_r|^2 + C \Big( \frac{1}{\ga^4}  \intl_{\Og \setminus B_\tau}  |v|^2  + \intl_{\Og \setminus B_\tau} |\mq|^2 \Big). 
\end{equation}
A combination of \eqref{E:Inf}, \eqref{HF-2} and \eqref{HF-5} yields
\begin{equation}\label{HF-6}
\frac{1}{2(d-1)} \intl_{B_{r_*} \setminus B_\tau} \Big( |\nabla v|^2 +\mf^2 \; |v|^2 \Big)   \le C \; \int\limits_{\Omega} |q|^2 +  \frac{r_* (3-d)}{2} \int\limits_{B_R \setminus B_{r_*}} \frac{u^2}{r^3} +  F_0(v) - F(R, v).
\end{equation}
Recall that   
\begin{equation*} - F(R,v) =  \; \Re \; \Big(\intl_{\partial B_R}  \frac{ {r_*}}{d-1} |v_r|^2  - \frac{{r_*}}{d-1}  |\nabla_{\partial B_R} v|^2 + \frac{{r_*}}{2R^2}  \; |v|^2 + \frac{{r_*}}{R} \, v_r \, \bar v  + \frac{\mf^2 \; {r_*}}{d-1} |v|^2 \Big).\end{equation*} 
Since $v$ satisfies the outgoing condition, it follows that \footnote{For the details of the argument, see the one used to obtain \cite[(2.19)]{NguyenVogelius1}.}
\begin{equation}\label{part1-1}
\limsup_{R \to \infty} \; -F(R,v) \leq \frac{2 \, \mf^2 \; {r_*}}{d-1} \; \limsup_{R \to \infty} \; \int\limits_{\pdh B_R}|v|^2.
\end{equation}
A combination of  \eqref{E:Inf} and \eqref{part1-1} yields  
\begin{equation} \label{E:FR} \limsup_{R \to \infty} \; -F(R,v) \leq  \frac{2 {r} \; \ga^2 \; \mf}{d-1} \; \int\limits_{\Omega} |q|^2 = C \eps^2 \; \int\limits_{\Omega} |q|^2 .\end{equation}
We have
\begin{equation} \label{E:F1}
F_0(v) = \Re \intl_{\partial B_\tau}  -\frac{2 \tau}{d-1} \; |\pn v|^2 + \frac{\tau}{d-1} \; |\nabla v|^2 + \pn v \; \overline v  - \frac{\mf^2 \; \tau}{d-1} \; |v|^2. 
\end{equation}
This implies
\begin{equation}  \label{E:F}
F_0(v) \leq  C \Big(\int\limits_{\partial B_{\tau}} |\nabla v|^2 +\mf^2 \int\limits_{\partial B_{\tau}} |v|^2 \Big).
\end{equation}
We derive from \eqref{HF-0} that 
\begin{equation}\label{E:FA} F_0(v) \; \leq\; C \; \int\limits_{\Omega} |q|^2.
\end{equation}
Combining (\ref{HF-6}), (\ref{E:FR}),  and (\ref{E:FA}), we obtain 
\begin{equation*}
\frac{1}{2(d-1)} \int\limits_{B_{r_*} \setminus B_\tau} |\nabla v|^2 +\mf^2 \; |v|^2 \leq C  \int\limits_{\Omega} |q|^2 +  \frac{r_*  (3-d)}{2} \int\limits_{B_R \setminus B_{r_*}} \frac{u^2}{r^3}. 
\end{equation*}
The proof for $d=3$ is complete. For $d=2$, it remains to absorb the second term on the RHS to the LHS. Without of generality we may assume that $r_*$ is big enough. Then, the absorption can be done as in \cite[pp. 11-12]{NguyenVogelius1}. The details are left to the reader. 
\proofend

\medskip

The following result will be used to obtain the estimate for $\mw_\ell$.

\begin{corollary}\label{cor-high2}
Let $d=2, 3$, $k_0> 0$, $r_0>0$, $0 < \eps < 1$, $\mf \geq k_0$, and $\mq \in L^2(\rN^d)$ with $\supp \mq \subset \overline \Og$.  Let $v \in H^1_{\loc}(\rN^d)$ be the unique outgoing solution to 
\begin{equation*} \Delta v+\mf^2\; v + i \,\mf \, \sigma_\eps \; v = \mq \mbox{ in } \rN^d.
\end{equation*} 
Then, for any $K \subset \subset \rN^d \setminus \bar\Og$, there is a positive constant $C_K$ independent of $\mf$, $\eps$, and $q$ such that
\begin{equation*}
\|\nabla v \|_{L^2(K)} + \mf  \|v\|_{L^2(K)}  \leq C_K \;\ga \;\mf^2 \;  \| q\|_{L^2(\Omega)} .
\end{equation*}
\end{corollary}
\noindent {\bf Proof.}
The proof of this corollary is similar to that of Corollary \ref{cor-high}. One only needs to use Lemma \ref{L:nhHF} (more precisely, the estimate \eqref{H:p3}) in place of Lemma \ref{lem-pro3-1}. The details are left to the reader.
\proofend

\subsection{Proof of Proposition~\ref{P:HF}} Proposition~\ref{P:HF} can now be proved in a similar way to Proposition~\ref{P:MF}.  We only need to use Corollary~\ref{cor-high} in place of Lemma~\ref{lem-stability1} and Corollary~\ref{cor-high2} in place of Lemma~\ref{lem-stability2}. We present the proof here for the convenience of the reader. 
From the definition of $\mw_\ell$ \eqref{def-md} and $\me_\ell$ \eqref{def-ed}, we have 
\begin{equation} \label{E:vweHF} \|v^\eps - v_{\ell}^{a}\|_{H^1(B_r \setminus \Og)} \leq \|\mw_\ell \|_{H^1(B_r \setminus \Og)} + \|\me_\ell \|_{H^1(B_r \setminus \Og)}. \end{equation} 
Applying Corollary~\ref{cor-high} and Lemma~\ref{L:ekMH} (with $m=1$), we have, for $\ell =0, 1$, 
\begin{equation} \label{E:melH} \|\nabla \me_\ell \|_{L^2(K)} + \mf \, \|\me_\ell \|_{L^2(K)} \leq C \,  \mf^2 \; \|h_\ell \|_{L^2(\rN^d)} \leq C \,\mf^{2\ell+3} \,  \ga^{\ell+1} \; \|\ms\|_{H^{2\ell+1}(\rN^d)}.\end{equation}
Applying Corollary~\ref{cor-high2} and Lemma~\ref{L:wkMH}, we obtain, for $\ell =0,  1, 2$, 
\begin{equation*}  \|\nabla \mw_\ell \|_{L^2(K)} + \mf \, \|\mw_\ell \|_{L^2(K)} \leq C\;\ga \;\mf^2 \;  \|q_\ell\|_{L^2(\Omega)} \leq C \; \mf^{2\ell+5} \; \ga^{\ell} \; \|\ms\|_{H^{2\ell+3}(\rN^d)}. \end{equation*}
This and Lemma~\ref{L:ueH} imply, for $\ell =0, 1$, 
\begin{align} \label{E:wl+3H}
\nonumber \|\nabla \mw_\ell \|_{L^2(K)} + \mf \, \|\mw_\ell \|_{L^2(K)}  \leq& \ga^{\ell+1} \, \big(\|\nabla w_e^\ell \|_{L^2(K)} + \mf \, \|w_e^\ell \|_{L^2(K)} \big) + C \,\mf^{2\ell+7} \,  \ga^{\ell + 1} \; \|\ms\|_{H^{2\ell+ 5}(\rN^d)} \\[6 pt] \leq& C \,\mf^{2\ell+7} \,  \ga^{\ell + 1} \; \|\ms\|_{H^{2\ell+ 5}(\rN^d)}.
\end{align}
A combination of \eqref{E:vweHF}, \eqref{E:melH}, and \eqref{E:wl+3H} yields the conclution. \proofend

\appendix
\section{Appendix: The Proof of the Estimates for the Asymptotic Expansions}
\renewcommand{\theequation}{A\arabic{equation}}
\renewcommand{\thelemma}{A\arabic{lemma}}
  \setcounter{lemma}{0}  

\subsection{Proof of Lemma \ref{L:ueM}}\label{S:ueM} 
The conclusion of Lemma~\ref{L:ueM} follows from the definition of $w_e^\ell$ for $\ell = 0, 1, 2$,  and the standard regularity theory of  elliptic equations. 
The details are left to the reader. \proofend

\subsection{Proof of Lemma~\ref{L:ueH}} \label{S:ueH}

Using the definition of $w_e^0$ and applying Lemma~\ref{L:ExtHF} for $r_*=R$, we have
\begin{equation*}
\|\nabla w_e^0 \|_{L^2(B_R \setminus \Og)} +\mf \; \|w_e^0\|_{L^2(B_R \setminus \Og)}  \leq  C \|s\|_{L^2(\rN^d  \setminus \Og)}. 
\end{equation*}
Using the regularity theory of elliptic equations and applying Lemma~\ref{L:ExtHF}, we have, for $m \ge 1$,  
\begin{equation*}
\|\nabla w_e^0 \|_{H^{m}(B_R \setminus \Og)} +\mf \; \|w_e^0\|_{H^m(B_R \setminus \Og)}  \leq  C k^{m} \|s\|_{H^m(\rN^d  \setminus \Og)}. 
\end{equation*}
and hence by the trace theory, we obtain
\begin{equation*}
\|\partial_n w_e^0\|_{H^{m-1/2}(\Gamma)} \le C  \|w_e^0\|_{H^{m+1}(\rN^d  \setminus \Og)} \leq  C k^{m+1} \|s\|_{H^{m+1}(\rN^d  \setminus \Og)}. 
\end{equation*}
By \eqref{def-wi1} and \eqref{def-wi2},  the conclusion follows from the definition of $w^1_e$, $w^2_e$ in \eqref{E:Recursive-e}, Lemma~\ref{L:ExtHF}, and the standard regularity  theory of elliptic equations. \proofend

\subsection{Proof of Lemma \ref{L:wkMH}}\label{S:wkMH} 
We only consider $\mf \geq k_0$. The other case follows similarly. From the definition of $q_{\ell}$, we have
\begin{equation*}  \mq_\ell(x) = \left\{\begin{array}{cl} \dsp - [\Delta v_i^\ell(x)+\mf^2\; v_i^\ell(x) + \frac{i}{\ga^2}\; v_i^\ell(x)]   - \ga^\ell [\Delta \varphi_{\ga}(x) +\mf^2\; \varphi_{\ga}(x) + \frac{i}{\ga^2}\; \varphi_{\ga}(x)] & \mbox{ in } \Og,  \\[6 pt]
0  & \mbox{ in }  \rN^d \setminus \Og.
\end{array}
\right.
\end{equation*}
It follows that
\begin{equation}\label{E:qlest}
\|\mq_\ell\|_{L^2(\rN^d)} \leq  \| \Delta v_i^\ell +\mf^2\; v_i^\ell + \frac{i}{\ga^2}\; v_i^\ell\|_{L^2(\Og)} +   \ga^\ell\Big\|\Delta \varphi_{\ga} +\mf^2 \varphi_{\ga} +  \frac{i}{\ga^2}\; \varphi_{\ga} \Big\|_{L^2(\Og)} 
\end{equation}
Since $\varphi_{\ga}(x) = \nu \; \chi(\nu/\ga) \; \pn w_e^\ell(x_\Ga)$, we obtain
\begin{multline*}  
\ga^\ell\Big\|\Delta \varphi +\mf^2 \varphi +  \frac{i}{\ga^2}\; \varphi \Big\|_{L^2(\Og)} \leq C \ga^{\ell} \Big(\|\pn w_e^\ell\|_{H^2(\Ga)}+ \mf^2  \; \|\pn w_e^\ell\|_{L^2(\Ga)}\Big)\\ \times  \,\big( \|\nu \chi(\nu/\ga)\|_{H^2(0,\infty)} + \frac{1}{\ga^2} \|\nu \chi(\nu/\ga)\|_{L^2(0,\infty)} \big).
\end{multline*}
We have \footnote{The scaling for the variable of $\chi$ gives us the optimal estimate in term of $\ga$.}
\begin{equation}  \label{E:optimal}
 \|\nu \chi(\nu/\ga)\|_{H^2(0,\infty)} + \frac{1}{\ga^2} \|\nu \chi(\nu/\ga)\|_{L^2(0,\infty)}  \leq C \ga^{-1/2}. 
\end{equation}
Therefore,
\begin{equation*}  
\ga^\ell\Big\|\Delta \varphi_{\ga} +\mf^2 \varphi_{\ga} +  \frac{i}{\ga^2}\; \varphi_{\ga} \Big\|_{L^2(\Og)} \leq C \, \ga^{\ell-1/2} \Big(\|\pn w_e^\ell\|_{H^2(\Ga)}+ \mf^2  \; \|\pn w_e^\ell\|_{L^2(\Ga)}\Big).
\end{equation*}
We derive from Lemma~\ref{L:ueH} with $m=3$ that, for $\ell = 0, 1, 2$, 
\ba \label{E:first}  
\ga^{\ell}\Big\|\Delta \varphi_{\ga} +\mf^2 \varphi_{\ga}   +  \frac{i}{\ga^2}\; \varphi_{\ga}  \Big\|_{L^2(\Omega)} \leq C\; \ga^{\ell-1/2} \mf^{2\ell+3} \; \|\ms\|_{H^{2\ell+3}(\rN^d)}.
\ea
Using \eqref{IMP} and the definitions of $v_{i}^1$ and $v_i^2$, as in  \cite{HJNg1}, we have
\begin{equation}\label{part1} \Big\| \Delta v_i^1+\mf^2 v_i^1 + \frac{i}{\ga^2}v_i^1 \Big\|_{L^2(\Og)} \leq \\C \left(\| {\cal A}_1 w_i^1 \|_{L^2(\Ga \times \rN_+)} + \ga^{-1}\|v_i^{1}\|_{H^1(\Ga \times \mR_+)}\right)
\end{equation}
and 
\begin{equation}\label{part2} \Big\| \Delta v_i^2+\mf^2 v_i^2 + \frac{i}{\ga^2}v_i^2 \Big\|_{L^2(\Og)} \leq \\C \left( \ga \| {\cal A}_1 w_i^2 +  {\cal A}_2 w_i^1 + {\cal A}_1 w_i^2 \|_{L^2(\Ga \times \rN_+)} + \|v_i^{2}\|_{H^1(\Ga \times \mR_+)}\right). 
\end{equation}
Recall that 
\begin{align*}
\mA_1 =& 2\mH  \partial_\eta + 6 \eta \mH  (\partial_\eta^2 + i), \\[6pt]
\mA_2 =& \Delta_\Gamma+\mf^2 + 2 \eta (\mG+ 4\mH^2)  \partial_\eta + 3 \eta^2 (\mG+ 4\mH^2) (\partial_\eta^2 +i). 
\end{align*}
Applying Lemma~\ref{L:ueH} with ($m =2$ and $\ell =0$) and using  \eqref{def-wi1}, we obtain:
\begin{equation}\label{E:vil1} 
\|w_i^1 \|_{L^2(\Ga \times \rN_+)} + \ga^{-1} \, \|v_i^{1}\|_{H^1(\Ga \times \mR_+)} \le C k^{2} \|s \|_{H^{2}(\rN^d)}. 
\end{equation}
Applying Lemma~\ref{L:ueH} with  ($m=3$ and $\ell = 0$) and ($m=1$ and $\ell =1$) and using  \eqref{def-wi2}, we obtain 

\begin{equation} \label{E:vil2} 
\|w_i^1 \|_{H^2(\Ga \times \rN_+)} + \|w_i^2 \|_{L^2(\Ga \times \rN_+)} \le C k^3 \|s \|_{H^{3}(\rN^d)}. 
\end{equation}
Moreover, applying Lemma~\ref{L:ueH} with  ($m=2$ and $\ell =0, 1$) 
\begin{equation} \label{E:vil3} 
\ga^{-1} \; \|v_i^{2}\|_{H^1(\Ga \times \mR_+)} =  \|w_i^1 + \ga \, w_i^2\|_{H^1(\Ga \times \mR_+)} \le C k^4 \|s \|_{H^{4}(\rN^d)}.
\end{equation} 
The conclusion now follows from \eqref{E:qlest}, \eqref{E:first}, \eqref{part1}, \eqref{part2}, \eqref{E:vil1}, \eqref{E:vil2}, \eqref{E:vil3}.  We note here that for the case $\ell=0$, we use the fact that $v_i^0 \equiv 0$. \proofend

\subsection{Proof of Lemma \ref{L:ekMH}} \label{S:ekMH}
 A calculation (see e.g., \cite[(4.34)]{HJNg1}) shows that $h_0 =0$, $h_1 = \frac{\ga^2}{\alpha} \; \pn w_e^1$. The conclusion follows from Lemmas~\ref{L:ueM} and \ref{L:ueH}.  \proofend

\bigskip
\noindent{\bf Acknowledgement.}  The first author thanks Patrick Joly and Houssem Haddar for introducing the subject. The first author was partially   
supported by NSF grant DMS-1201370 and by the Alfred P. Sloan Foundation. The second author was partially supported by NSF grant DMS-1212125.



\end{document}